\documentclass[a4paper]{amsart}

\usepackage{amsfonts}
\usepackage{bbm}
\usepackage{amssymb}
\usepackage{amsthm}
\usepackage{amsmath}
\usepackage{amscd}
\usepackage{graphicx}
\usepackage{stmaryrd}
\usepackage{scrpage2}
\usepackage[all]{xy}
\usepackage{tabularx}
\SelectTips{cm}{}
\DeclareMathOperator{\NIL}{NIL}
\DeclareMathOperator{\Nil}{Nil^{\fd}}
\DeclareMathOperator{\Mat}{M}
\DeclareMathOperator{\GL}{GL}
\DeclareMathOperator{\colim}{colim}
\DeclareMathOperator{\NA}{NA^{\fd}_+}
\DeclareMathOperator{\NK}{NK}
\DeclareMathOperator{\NAM}{NA_{--}^{\fd}}
\DeclareMathOperator{\NAPM}{NA^{\fd}_\pm}
\DeclareMathOperator{\fd}{fd}
\DeclareMathOperator{\f}{f}
\DeclareMathOperator{\TOP}{TOP}
\newlength{\wlength}
\newtheorem{Th}{Theorem}[section]
\newtheorem*{Th*}{Theorem}
\newtheorem{Pro}[Th]{Proposition}
\newtheorem*{Pro*}{Proposition}
\newtheorem{Le}[Th]{Lemma}
\newtheorem{Co}[Th]{Corollary}
\newtheorem*{Co*}{Corollary}

\theoremstyle{definition}

\newtheorem{De}[Th]{Definition}
\newtheorem*{De*}{Definition}
\newtheorem{Rem}[Th]{Remark}

\newtheorem{No}[Th]{Notation}

\newcommand{\al}{\alpha}
\newcommand{\zz}{\mathbb{Z}}

\newcommand{\nn}{\mathbb{N}}

\newcommand{\xym}{\xymatrix}
\newcommand{\ra}{\rightarrow}

\newcommand{\lra}{\longrightarrow}
\newcommand{\ot}{\otimes}

\renewcommand{\SS}{\mathbb{S}}
\renewcommand{\AA}{\mathbb{A}}
\newcommand{\BB}{\mathbb{B}}
\newcommand{\CC}{\mathbb{C}}
\newcommand{\DD}{\mathbb{D}}

\newcommand{\HH}{\mathbb{H}}

\newcommand{\sS}{\mathcal{S}}
\newcommand{\sC}{\mathcal{C}}

\newcommand{\NN}{\mathbb{N}}
\renewcommand{\SS}{\mathbb{S}}
\newcommand{\ZZ}{\mathbb{Z}}
\newcommand{\FF}{\mathbb{F}}

\newcommand{\TC}{\mathrm{TC}}
\newcommand{\THH}{\mathrm{THH}}
\newcommand{\Wh}{\mathrm{Wh}}
\renewcommand{\GL}{\mathrm{GL}}
\newcommand{\id}{\mathrm{id}}
\newcommand{\Map}{\mathrm{Map}}
\newcommand{\Nat}{\mathrm{Nat}}
\newcommand{\ev}{\mathrm{ev}}
\newcommand{\trc}{\mathrm{trc}}

\newcommand{\hofiber}{\mathrm{hofiber}}

\newcommand{\pco}{^{\wedge}_p}

\begin{document}

\title{Operations on the A-theoretic nil-terms}

\author{Joachim Grunewald, John R. Klein and Tibor Macko}

\date{\today}



\thanks{J.G and T.M. are supproted by SFB 478 Geometrische Strukturen in der Mathematik, M\"unster.}

\address{Mathematisches Institut \\ Universit\"at M\"unster \\
Einsteinstrasse 62 \\ M\"unster, D-48149 \\ Germany}
\email{grunewal@math.uni-muenster.de}
\address{Wayne State University \\ Detroit \\ MI 48202 \\ USA}
\email{klein@math.wayne.edu}
\address{Mathematisches Institut \\ Universit\"at M\"unster \\
Einsteinstrasse 62 \\ M\"unster, D-48149 \\ Germany \\ and
Matematick\'y \'Ustav SAV \\ \v Stef\'anikova 49 \\ Bratislava,
SK-81473 \\ Slovakia} \email{macko@math.uni-muenster.de}

\numberwithin{equation}{section}


%
\maketitle

\begin{abstract}
For a space $X$, we define Frobenius and Verschiebung operations on
the nil-terms $\NAPM (X)$ in the algebraic $K$-theory of spaces, in
three different ways. Two applications are included. Firstly, we
show that the homotopy groups of $\NAPM (X)$ are either trivial or
not finitely generated as abelian groups. Secondly, the Verschiebung
defines a $\ZZ[\NN_\times]$-module structure on the homotopy groups
of $\NAPM (X)$, with $\NN_\times$ the multiplicative monoid.

We also give a calculation of the homotopy type of the nil-terms
$\NAPM (\ast)$ after $p$-completion for an odd prime $p$ and their
homotopy groups as $\ZZ_p [\NN_\times]$-modules up to dimension
$4p-7$. We obtain non-trivial groups only in dimension $2p-2$, where
it is finitely generated as a $\ZZ_p[\NN_\times]$-module, and in
dimension $2p-1$, where it is not finitely generated as a $\ZZ_p
[\NN_\times]$-module.
\end{abstract}

\section{Introduction}
The fundamental theorem of algebraic $K$-theory of spaces which is
proved in \cite{fund1} states that for any space $X$ there is a
splitting
\begin{equation} \label{A-fund-thm}
A^{\fd} (X \times S^1) \simeq A^{\fd}(X) \times B A^{\fd}(X) \times
\NA (X) \times \NAM (X),
\end{equation}
where $A^{\fd}(X)$ denotes the finitely-dominated version of the
algebraic $K$-theory of the space $X$, by $B A^{\fd} (X)$ is denoted
certain canonical non-connective delooping of $A^{\fd} (X)$ and
$\NA(X)$, $\NAM(X)$ are two homeomorphic nil-terms. Thus the study
of $A^{\fd}(X \times S^1)$ splits naturally into studying
$A^{\fd}(X)$ and the nil-terms. Over time there has been steady
progress in understanding $A^{\fd}(X)$ for some $X$, mostly for $X =
\ast$ \cite{Klein-Rognes(1997)}, \cite{Rognes(2002)},
\cite{Rognes(2003)}, but not much is known about the nil-terms.
These are the subject of the present paper.

Recall that the splitting (\ref{A-fund-thm}) is analogous to the
splitting of the fundamental theorem of algebraic $K$-theory of
rings \cite{Grayson(1976)} for any ring $R$
\begin{equation} \label{K-fund-thm}
K(R[\ZZ]) \simeq K(R) \times B K(R) \times \NK_+ (R) \times \NK_-
(R),
\end{equation}
where $K(R)$ denotes the $(-1)$-connective $K$-theory space of the
ring $R$ and $B K(R)$ is a certain canonical non-connective
delooping of $K(R)$. In fact, for any space $X$ there is a {\it
linearization map} $l \colon A^{\fd} (X) \ra K(\ZZ[\pi_1 X])$ and
the two splittings are natural with respect to this map.

We pursue the analogy between the algebraic $K$-theory of spaces
(the non-linear situation) and the algebraic $K$-theory of rings
(the linear situation) further. We define the Frobenius and
Verschiebung operations on the nil-terms $\NAPM (X)$, which are
analogs of such operations defined in the linear case. As a
consequence we obtain the following result (compare with
\cite{Farrell}):

\begin{Th} \label{thm-1} \label{infinite}
The homotopy groups $\pi_\ast \NAPM (X)$ and all of their
$p$-primary subgroups are either trivial or not finitely generated
as abelian groups.
\end{Th}

This result suggests that it is very hard to determine the homotopy
groups of $\NAPM (X)$ and that it is important to have more
structure on them. We achieve this in Corollary \ref{ZZ[NN]-module}
which says that the Verschiebung operations define a structure of a
$\ZZ[\NN_\times]$-module on the homotopy groups of the nil-terms.
Here $\NN_\times = \{1, 2, \ldots \}$ is a monoid with respect to
multiplication. In the algebraic $K$-theory of rings the
corresponding structure was studied in \cite{Conolly-et-al(1995)},
\cite{Loday(1980)}.

It is standard knowledge that $\NK_\pm (\ZZ) \simeq \ast$, since
$\ZZ$ is a regular ring \cite{Rosenberg(1994)}. However, not much is
known about the spaces $\NAPM (\ast)$. One approach is to compare
them with the linear nil-terms via the linearization map $\NAPM
(\ast) \ra \NK_\pm (\ZZ)$. This is known to be a rational
equivalence. Moreover, after $p$-completion at a prime $p$ it is
$(2p-3)$-connected. Hence the nil-terms $\NAPM (\ast)$ are
rationally trivial. It seems that the only other known results in
this direction are that they are trivial in dimensions $0$, $1$ and
non-trivial in dimension $2$ \cite{Waldhausen(1978)},
\cite{Igusa(1982)}.

We give a calculation of the homotopy type of the nil-terms $\NAPM
(\ast)$ and of the module structure on their homotopy groups in a
certain range. For this we first $p$-complete $\NAPM (\ast)$ at an
odd prime $p$ so that the homotopy groups become
$\ZZ_p[\NN_\times]$-modules. Below $\HH \FF_p$ denotes the
Eilenberg-Mac Lane spectrum of $\FF_p$ and $S^1 (n)$ is just a copy
of $S^1$ indexed by $n \in \NN_\times$. The result is:

\begin{Th} \label{pi-Nil-calculation}
If $p$ is an odd prime the following holds.
\begin{enumerate}
\item There is a $(4p-7)$-connected map
\[
\bigvee_{\pm n \in \NN_\times} \Sigma^{2p-2} \HH \FF_p \wedge
(S^1(\pm n)_+) \lra \NAPM (\ast)\pco~,
\]
hence
\begin{align*}
\pi_{2p-2} \NAPM (\ast)\pco & \cong \oplus_{n \in \NN_\times} \FF_p \{\beta_{\pm n}\} \\
\pi_{2p-1} \NAPM (\ast)\pco & \cong \oplus_{n \in \NN_\times} \FF_p \{\gamma_{\pm n}\} \\
\pi_i \NAPM (\ast)\pco & \cong 0 \; \textup{for } i < 2p-2 \textup{
, } 2p-1 < i \leq 4p-7,
\end{align*}
where
\begin{align*}
\beta_{\pm n} & \in \pi_{2p-2} \Sigma^{2p-2} \HH \FF_p \wedge
(S^1(\pm n)_+) \cong \FF_p \\
\gamma_{\pm n} & \in \pi_{2p-1} \Sigma^{2p-2} \HH \FF_p \wedge
(S^1(\pm n)_+) \cong \FF_p
\end{align*}
represent a certain choice of generators of these $\ZZ_p$-modules.
\item The $\ZZ_p[\NN_\times]$-module structure on $\pi_\ast \NAPM (\ast)\pco$
is given by
\begin{align*}
(n , \beta_m) & \mapsto \beta_{nm} \\
(n , \gamma_m) & \mapsto n \cdot \gamma_{nm}
\end{align*}
\end{enumerate}
\end{Th}

The proof of Theorem \ref{pi-Nil-calculation} follows the general
strategy suggested by Madsen in \cite{Madsen(1995)} to study the
linearization map $A^{\fd} (S^1) \ra K(\ZZ[\ZZ])$. Certain homotopy
theoretic constructions and relations between $A^{\fd} (X)$, the
trace invariants of $K$-theory, and the free loop space of $X$ are
used.

In the linear situation there are partial results pointing in the
direction that the homotopy groups of the nil-terms could be
finitely-generated as $\ZZ[\NN_\times]$-modules, see
\cite{Conolly-et-al(1995)}, \cite{Loday(1980)}. However, Theorem
\ref{pi-Nil-calculation} shows that this is not the case in the
non-linear situation. We have

\begin{Co} \label{ZZ[NN]-module-nil}
The $\ZZ_p[\NN_\times]$-module $\pi_{2p-2} \NAPM (\ast)\pco$ is
(finitely) generated by $\beta_{\pm 1}$, and the
$\ZZ_p[\NN_\times]$-module $\pi_{2p-1} \NAPM (\ast)\pco$ is not
finitely generated.
\end{Co}

Combining the Frobenius and Verschiebung operations we can define
also a structure of a $\ZZ_p [\NN_\times \ast \NN_\times]$-module on
the homotopy groups of $\NAPM(X)\pco$, where $\NN_\times \ast
\NN_\times$ denotes the free product of monoids. In Corollary
\ref{ZZ[NN*NN]-module-F_K} we show that $\pi_{2p-1} \NAPM
(\ast)\pco$ is also not finitely-generated as a $\ZZ_p [\NN_\times
\ast \NN_\times]$-module.

For a geometric application of our results recall that algebraic
$K$-theory of spaces is related via the assembly maps to
automorphisms of manifolds. For a general scheme we refer the reader
to the survey article \cite{Weiss-Williams(2001)}. This general
scheme can be applied for example in the case of smooth manifolds
with negative sectional curvature. For such a manifold $M$ of
dimension $n$ the work of Farrell and Jones
\cite{Farrell-Jones(1991)} gives the following results about the
space $\TOP (M)$ of self-homeomorphisms of $M$:
\[
\pi_j \TOP (M) \cong \bigoplus_T \pi_{j+2} \NA(\ast)
\]
for $1 < j \leq \phi (n)$, where $T$ runs through simple closed
geodesics in $M$ and $\phi(n)$ denotes the concordance stable range
which is approximately $n/3$, and there is an exact sequence
\begin{multline*}
0 \ra \bigoplus_T \pi_3 \NA (\ast) \ra \pi_1 \TOP (M) \ra Z(\pi_1 M) \ra \\
\ra \bigoplus_T \pi_2 \NA (\ast) \ra \pi_0 \TOP (M) \ra \mathrm{Out}
(\pi_1 M).
\end{multline*}

So our results describe new non-trivial families of automorphisms of
negatively curved smooth manifolds.


\noindent\textbf{Organization.} Section 2 contains the background
about algebraic $K$-theory and nil-terms. In section 3 we define the
Frobenius and Verschiebung operations in three different ways. In
section 4 we prove that the three definitions coincide. In section 5
we prove the certain identities satisfied by these operations and
derive some consequences out of them. Section 6 contains the proof
of Theorem 1.2.

\noindent\textbf{Acknowledgements.} We would like to thank John
Rognes for help with the earlier version of this paper and for some
arguments in section 6 of the present version. Further we have
benefited from conversations with Tilman Bauer, Andr\'e Henriques,
Steffen Sagave and Marco Varisco. The second author is indebted to
Tom Goodwillie for discussions in connection with the material in
section 6.


\section{Preliminaries}

Waldhausen defined algebraic $K$-theory for categories with
cofibrations and weak equivalences \cite{Waldhausen1126}. For a
topological space $X$ its algebraic $K$-theory $A (X)$ can be
defined as the algebraic $K$-theory of the category of based spaces
with an action of the geometric realization of the Kan loop group of
$X$ or using equivariant stable homotopy theory, in terms of certain
spaces of matrices. The modern approach to the latter alternative
uses the language of $\SS$-algebras and is presented in
\cite{Dundas-Goodwillie-McCarthy(2004)}. We briefly review the
definitions in both cases. Further we recall two definitions of the
nil-terms which are the principal objects of study in this paper.

\subsection{K-theory}

We follow the notation of \cite{fund1}, \cite{fund2}, \cite{fund3}.
Recall that a {\it Waldhausen category} is a category with
cofibrations and weak equivalences satisfying axioms of
\cite[1.2]{Waldhausen1126}. If $\sC$ is such a category there is
defined a connected based space $|w\sS_\bullet \sC|$ (in fact an
infinite loop space), the $\sS_\bullet$-construction of $\sC$ whose
loop space is taken to be the definition of the algebraic $K$-theory
of $\sC$.

We recall the examples of Waldhausen categories relevant for us. Let
$M_\bullet$ be a simplicial monoid whose realization is denoted by
$M$. Let $\mathbb{T}(M)$ be the category whose objects are based
$M$-spaces, i.e. based spaces $Y$ equipped with a based left action
$M \times Y  \rightarrow Y$. Morphisms of $\mathbb{T}(M)$ are based
equivariant maps. The cell of dimension $n$ is defined to be
\[
D^n \times M.
\]
If $Z$ is an object in $\mathbb{T}(M)$, and $\al \colon S^{n-1}
\times M \rightarrow Z$ is an equivariant map we can form the
pushout
\[
Z \cup_\al (D^n \times M).
\]
This operation is called the effect of \emph{attaching a cell}. A
morphism $Y \rightarrow Z$ of $\mathbb{T}(M)$ is a
\emph{cofibration} if either $Z$ is obtained from $Y$ by a sequence
of cell attachments or it is a retract of the forgoing. An object is
said to be \emph{cofibrant} if the inclusion of the base point is a
cofibration.

Let $\mathbb{C}(M) \subset \mathbb{T}(M)$ be the full subcategory of
cofibrant objects. An object of $\mathbb{C}(M)$ is called
\emph{finite} if it is isomorphic to a finite free based
$M$-CW-complex, i.e. it is built from a point by a finite number of
cell attachments. An object $Y$ is called \emph{homotopy finite} if
there is a finite object $X$ and a morphism $Y \rightarrow X$ which
considered as a map of ordinary spaces is a weak homotopy
equivalence. An object is called \emph{finitely dominated} if it is
a retract of a homotopy finite object. The full subcategory of
finite (finitely dominated) objects, is denoted $\mathbb{C}_{\f}
(M)$ ($\mathbb{C}_{\fd}(M)$). Both categories $\mathbb{C}_{\f} (M)$
and $\mathbb{C}_{\fd}(M)$ have the structure of a Waldhausen
category by defining a morphism to be a cofibration if it is a
cofibration in $\mathbb{T}(M)$ and a weak homotopy equivalence if it
is a weak homotopy equivalence of the underlying spaces.

\begin{No}
Denote
\begin{align*}
A^{\fd}(*,M) & := \Omega |h \mathcal{S}_\bullet
\mathbb{C}_{\fd}(M)|,
\\
A^{\f}(*,M) & := \Omega |h \mathcal{S}_\bullet \mathbb{C}_{\f} (M)|.
\end{align*}
\end{No}

\begin{Rem}
Recall that if $X$ is a connected and based topological space and
$M$ is the geometric realization of its Kan loop group, then $A^{\f}
(\ast,M)$ is one of the definitions of $A(X)$ given in
\cite{Waldhausen1126}. In analogy we regard $A^{\fd} (\ast,M)$ as
the definition of $A^{\fd} (X)$.
\end{Rem}

\begin{Rem}
The relationship between the two versions is
\begin{equation*}
A^{\fd}(*,M) \simeq \tilde{K}_0(\zz [ \pi_0 M ]) \times A^{\f}(*,M),
\end{equation*}
where $\tilde{K}_0(\zz [ \pi_0 M ])$ is the reduced class group of
the monoid ring $\ZZ [\pi_0 M]$, \cite[Lemma 1.7.(3)]{fund1}.  It is
necessary to use the finitely dominated version in order for the
fundamental theorem to hold.
\end{Rem}

For the second definition of $A^{\fd} (X)$ we need the language of
$\SS$-algebras, which we briefly recall following
\cite{Dundas-Goodwillie-McCarthy(2004)}. By $\Gamma^{op}$ is denoted
the category of finite pointed sets and pointed maps, more precisely
the skeletal subcategory consisting of the sets $k_+ = \{0, \ldots ,
k\}$, where $0$ is the base point. By $\sS_\ast$ is denoted the
category of pointed simplicial sets. A $\Gamma${\it -space} $\BB$ is
a functor $\BB \colon \Gamma^{op} \ra \sS_\ast$. The category of
$\Gamma$-spaces is denoted $\Gamma\sS_\ast$, the morphisms are
natural transformations of functors from $\Gamma^{op}$ to
$\sS_\ast$. $\Gamma$-spaces give rise to spectra, cf.
\cite[II.2.1.12]{Dundas-Goodwillie-McCarthy(2004)}. There are
several relevant examples for us. The {\it sphere spectrum} is $\SS
\colon k_+ \mapsto k_+$, where $k_+$ in the target is a constant
pointed simplicial set. For any $X \in \sS_\ast$ the {\it suspension
spectrum} is given by $\Sigma^\infty X \colon k_+ \mapsto X \wedge
k_+$ and for an abelian group $G$ the {\it Eilenberg-MacLane
spectrum} $\HH G$ is given by $k_+ \mapsto G^{\times k}$, cf.
\cite[II.1.2]{Dundas-Goodwillie-McCarthy(2004)}. The homotopy groups
of a $\Gamma$-space are defined as the homotopy groups of the
corresponding spectrum.

There is a notion of a smash product of $\Gamma$-spaces $\BB \wedge
\CC$, cf \cite[II.1.2.3]{Dundas-Goodwillie-McCarthy(2004)},
\cite{Lydakis(1999)}, and a unit of the smash product, which is the
already mentioned sphere spectrum $\SS$. The triple $(\Gamma
\sS_\ast,\wedge, \SS)$ turns out to be a symmetric monoidal
category. An $\SS${\it -algebra} $\AA$ is a monoid in $(\Gamma
\sS_\ast,\wedge,\SS)$, i.e. a $\Gamma$-space with a multiplication
map $\mu \colon \AA \wedge \AA \ra \AA$ and a unit map $1 \colon \SS
\ra \AA$ satisfying certain relations, cf
\cite[II.1.4.1]{Dundas-Goodwillie-McCarthy(2004)}. $\SS$-algebras
give rise to ring spectra with a well-behaved smash product. The
relevant examples here start again with the sphere spectrum $\SS$
with the multiplication map given by the action of $\SS$ on itself.
For a simplicial monoid $M$ there is the {\it monoid }$\SS${\it
-algebra} $\SS[M] = \Sigma^\infty(M_+) \colon k_+ \ra M_+ \wedge
k_+$ with the multiplication map induced by the multiplication in
$\SS$ and the monoid operation of $M$. For a ring $R$ there is the
Eilenberg-MacLane $\SS$-algebra $\HH R$, with the multiplication
given by the multiplication in $R$. The homotopy groups of an
$\SS$-algebra form a simplicial ring.

For an $\SS$-algebra $\AA$ the $\SS$-algebra $\Mat_m(\AA)$ of $m
\times m$ matrices is
\begin{equation*}
\Mat_m (\AA) \colon k_+ \mapsto \Map (m_+, m_+ \wedge \AA(k_+)) =
\prod_0^{m-1} \bigvee_0^{m-1} \AA(k_+),
\end{equation*}
with the multiplication $\Mat_m (\AA) \wedge \Mat_m (\AA) \ra \Mat_m
(\AA)$ given in \cite[Example
1.4.4.(6)]{Dundas-Goodwillie-McCarthy(2004)}. The ``points" in
$\Mat_m(\AA)(k_+)$ are matrices with entries in $\AA(k_+)$ with at
most one non-zero entry in every column. By $\widehat{\Mat}_m (\AA)$
is denoted the simplicial ring $\Omega^\infty \Mat_m (\AA)$ whose
group of components is $\pi_0 \widehat{\Mat}_m(\AA) = \Mat_m (\pi_0
(\AA))$. The simplicial ring of ``invertible matrices"
$\widehat{\GL}_m(\AA)$ is defined as the pullback
\[
\xymatrix{ \widehat{\GL}_m (\AA) \ar[r] \ar[d] & \widehat{\Mat}_m(\AA) \ar[d] \\
\GL_m(\pi_0 (\AA)) \ar[r] & \Mat_m (\pi_0 (\AA)).}
\]

\begin{De}{\cite[III.2.3.2]{Dundas-Goodwillie-McCarthy(2004)}}
The $K$-theory of an $\SS$-algebra $\AA$ is defined as
\begin{equation*}
K(\AA) = K_0 (\ZZ [\pi_0 \AA]) \times \colim_{m \ra \infty} B
\widehat{\GL}_m (\AA)^+,
\end{equation*}
where the superscript $^+$ denotes the Quillen $+$-construction.
\end{De}
It can be shown that $K(\AA)$ is an infinite loop space, so we think
of $\AA \mapsto K(\AA)$ as of a functor from $\SS$-algebras to
spectra.

\begin{Rem} \label{ident-A-K}
The examples of the $\SS$-algebras mentioned above yield for a based
space $X$ with the Kan loop group $\Omega X$ and for a ring $R$
\begin{align*}
K (\SS[\Omega X]) & \simeq A^{\fd} (X) = A^{\fd}(\ast,\Omega X), \\
K (\HH R) & \simeq K(R),
\end{align*}
where $K(R)$ denotes the Quillen $K$-theory of a ring $R$ (more
precisely its (-1)-connective version), cf \cite[Chapter
III]{Dundas-Goodwillie-McCarthy(2004)}.
\end{Rem}

Finally we recall the {\it linearization map}. This is for any space
$X$ with the Kan loop group $\Omega X$ the map of $\SS$-algebras $l
\colon \SS[\Omega X] \ra \HH \ZZ [\pi_1 X]$ given by sending a loop
in $\Omega X$ to its equivalence class in $\pi_1 X$. This map
induces the map $l \colon A^{\fd} (X) \ra K(\ZZ [\pi_1 X])$ which we
also call the linearization map.

\subsection{Nil-terms}

We closely follow \cite{fund1}. For the definition of the nil-terms
we need to recall the definition of the mapping telescope from
\cite[page 27]{fund1}. The symbols $\nn_+$, $\nn_-$ or $\zz$ will
always denote monoids with the addition and with the generators $t$,
$t^{-1}$, or $t$. For an object $Y_+ \in \mathbb{C}_{\fd}(M \times
\nn_+)$ the object $Y_+(t^{-1}) \in \mathbb{C}_{\fd} (M \times \zz)$
is defined to be the categorical colimit of the sequence
\[
\xym{\cdots \ar[r]^-{t} & Y_+ \ar[r]^-{t} & Y_+ \ar[r]^-{t} & \cdots
.}
\]
Explicitly it is the quotient space of $Y_+ \times \zz$ in which a
pair $(y,n)$ is identified with the pair $(t \cdot y,n+1)$. The
action of $t^k$ for $k \in \zz$ on a pair $(y,n)$ yields the pair
$(y,n-k)$. This assignment defines an exact functor \cite[page
26]{fund1}.

Similarly, for $Y_- \in \mathbb{C}_{\fd}(M \times \nn_-)$, we have
the mapping telescope $Y_-(t) \in \mathbb{C}_{\fd}(M \times \zz)$
given by the colimit of the sequence
\[
\xym{\cdots \ar[r]^-{t^{-1}} & Y_- \ar[r]^-{t^{-1}} & Y_-
\ar[r]^-{t^{-1}} & \cdots .}
\]
As above we use the quotient of $Y_- \times \zz$ in which a pair
$(y,n)$ is identified with the pair $(t^{-1} \cdot y,n+1)$. The
action of $t^k$ for $k \in \zz$ on a pair $(y,n)$ yields the pair
$(y,n-k)$.

Recall from \cite[page 36]{fund1} that $\mathbb{D}_{\fd}(M \times
L)$ is defined to be the category whose objects are diagrams
\[
\xym{Y_- \ar[r] & Y & Y_+ \ar[l]}
\]
with $Y_- \in \mathbb{C}_{\fd}(M \times \nn_-)$, $Y \in
\mathbb{C}_{\fd}(M \times \zz)$ and $Y_+ \in \mathbb{C}_{\fd}(M
\times \nn_+)$, and where the maps $Y_- \rightarrow Y$ and $Y_+
\rightarrow Y$ are assumed to be cofibrations. We take the liberty
of specifying an object as a diagram or as a triple $(Y_-,Y,Y_+)$.

A morphism $(Y_-,Y,Y_+) \rightarrow (Z_-,Z,Z_+)$ of
$\mathbb{D}_{\fd}(M \times \zz)$ is a morphism $Y_- \rightarrow
Z_-$, a morphism $Y \rightarrow Z$ and a morphism $Y_+ \rightarrow
Z_+$ so that the evident diagram commutes. A cofibration in
$\mathbb{D}_{\fd}(M \times \zz)$ is a morphism such that each of the
maps $Y_- \rightarrow Z_-$, $Y_+ \rightarrow Z_+$ and $Y \rightarrow
Z$ is a cofibration and the induced maps
\begin{align*}
Y \cup_{Y_-(t)} Z_-(t) & \rightarrow Z \\
Y \cup_{Y_+(t^{-1})} Z_+(t^{-1}) & \rightarrow Z
\end{align*}
are cofibrations.

Further $\mathbb{D}_{\fd}(M \times \nn_\pm) \subset
\mathbb{D}_{\fd}(M \times \zz)$ is defined to be the full
subcategory consisting of objects with the additional property that
the induced map
\[
Y_\mp (t) \rightarrow Y
\]
is a weak equivalence.

For $L=\nn_+$, $\nn_-$ or $\ZZ$ recall from \cite[page 48]{fund1} an
augmentation functor
\begin{align*}
\epsilon \colon \mathbb{D}_{\fd}(M \times L) & \rightarrow \mathbb{C}_{\fd}(M) \\
(Y_-,Y,Y_+) & \mapsto Y/\zz.
\end{align*}
The homotopy fiber of the induced map
\[
\epsilon \colon  |h \mathcal{S}_\bullet \mathbb{D}_{\fd}(M \times
L)| \rightarrow |h \mathcal{S}_\bullet \mathbb{C}_{\fd}(M)|
\]
is denoted by $|h \mathcal{S}_\bullet \mathbb{D}_{\fd}(M \times
L)|^\epsilon$. The canonical map
\[
|h \mathcal{S}_\bullet \mathbb{D}_{\fd}(M \times L)|^\epsilon
\rightarrow |h \mathcal{S}_\bullet \mathbb{D}_{\fd}(M \times L)|
\]
admits a left homotopy inverse \cite[page 48]{fund1} and there is a
homotopy equivalence
\[
|h \mathcal{S}_\bullet \mathbb{D}_{\fd}(M \times L)| \simeq |h
\mathcal{S}_\bullet \mathbb{D}_{\fd}(M \times L)|^\epsilon \times |h
\mathcal{S}_\bullet \mathbb{C}_{\fd}(M)|.
\]
We have
\begin{De}{\cite[page 51]{fund1}}
\begin{align*}
\NA(*,M)  : & = \Omega |h \mathcal{S}_\bullet \mathbb{D}_{\fd} (M \times \nn_+)|^\epsilon \\
\NAM(*,M) : & = \Omega |h \mathcal{S}_\bullet \mathbb{D}_{\fd} (M
\times \nn_-)|^\epsilon.
\end{align*}
\end{De}

Similarly as in the linear case the nil-terms can be identified with
a subgroup of the $K$-theory of a certain nil-category. Here we
follow \cite{fund3}. As before let $M$ be the geometric realization
of a simplicial monoid $M_\bullet$. The nil-category
$\NIL_{\fd}(*,M)$ has as its objects pairs $(Y,f)$ where $Y$ is an
object in $\mathbb{C}_{\fd}(M)$ and $f$ is an $M$-map from $Y$ to
$Y$ which is homotopically nilpotent, i.e. there exist a non
negative integer $k$ such that $f^{\circ k}$ is equivariantly null
homotopic. A morphism $(Y,f) \rightarrow (Y',f')$ is a map $e \colon
Y \rightarrow Y'$ in $\mathbb{C}_{\fd}(M)$ such that $f' \circ e = e
\circ f$. A \emph{cofibration} in $\NIL_{\fd}(*,M)$ is a morphism
whose underlying map in $\mathbb{C}_{\fd}(M)$ is a cofibration. A
\emph{weak homotopy equivalence} is a morphism whose underlying map
in $\mathbb{C}_{\fd}(M)$ is a weak equivalence. This gives
$\NIL_{\fd}(*,M)$ the structure of a Waldhausen category \cite[Lemma
1.2]{fund3} and therefore it has a $K$-theory $K(\NIL_{\fd}
(\ast,M))$.
The forgetful functor
\begin{align*}
\NIL_{\fd}(*,M) & \rightarrow \mathbb{C}_{\fd}(M) \\
(Y,f) & \mapsto Y
\end{align*}
induces a map in $K$-theory whose homotopy fiber is denoted $\Nil
(\ast,M)$.
The main theorem of \cite{fund3} says
\begin{Th}[\cite{fund3}] \label{thm-fund3}
There are natural homotopy equivalences
\[
\xymatrix{ \delta_\pm \colon \Nil (*,M) \ar[r]^{\simeq} &  \Omega
\NAPM (*,M). }
\]
\end{Th}

\begin{Rem}
Another possibility is to define the nil-terms in the language of
$\SS$-algebras as the fibers of the map in $K$-theory induced by the
augmentation map of $\SS$-algebras $\SS[M \times L] \ra \SS[M]$
where $L = \NN_+$ or $\NN_-$. We will not use this definition in
this paper, so we omit the details.
\end{Rem}

\section{Frobenius and Verschiebung Operations}
In this section we define for a natural number $n \in \NN_\times$
the operations Frobenius $F_n$ and Verschiebung $V_n$ in the
algebraic $K$-theory of spaces using all three definitions from the
previous section. We also show that they restrict to operations on
the nil-terms.

Throughout the symbol $L$ will denote either $\nn_+$, $\nn_-$ or
$\zz$ as an additive monoid with the generator $t$, $t^{-1}$ or $t$.
For $n \in \NN_\times$ we will use the monoid homomorphism
\begin{align*}
\varphi_n \colon & L \ra L \\
& l \mapsto n \cdot l.
\end{align*}
Informally, in all three cases, the Verschiebung operation $V_n$
will be induced by $\varphi_n$ and the Frobenius operation $F_n$
will be the corresponding transfer.

\subsection{Operations on $A^{\fd}(*,M \times L)$}
\label{sec:op-on-na}
The definition is on the categorical level and goes in three steeps.
First we define the operations on $\mathbb{C}_{\fd}(M \times L)$,
then on $\mathbb{D}_{\fd}(M \times L)$ and then we show that the
operations restrict to the nil-terms.
\begin{De}[Verschiebung]
Define an exact functor
\begin{align*}
\text{V}_n \colon \mathbb{C}_{\fd}(M \times L) & \rightarrow
\mathbb{C}_{\fd}(M \times L), \\
Y & \mapsto Y \otimes_{\varphi_n} L,
\end{align*}
where $Y \otimes_{\varphi_n} L = Y \times L / \sim$ such that $(*,k)
\sim (*,l)$ and $(y,\varphi_n (k) + l) \sim (t^k \cdot y,l)$ for all
$k,l$. The $L$-action on $Y \ot_{\varphi_n} L$ is given by
$(t^k,(y,l)) \mapsto (y,k+l)$. The induced endomorphisms of
$A^{\fd}(*, M \times L)$ are also denoted by V$_n$ and are called
the \emph{Verschiebung operations}.
\end{De}
\begin{De}[Frobenius]
Define an exact functor
\begin{align*}
\text{F}_n \colon \mathbb{C}_{\fd}(M \times L) & \rightarrow
\mathbb{C}_{\fd}(M \times L) \\
Y & \mapsto Y_{\varphi_n}
\end{align*}
where $Y_{\varphi_n}$ is the same space $Y$ but with the $L$-action
precomposed with $\varphi_n$, i.e. $(t^k,y) \mapsto t^{\varphi_n
(k)} \cdot y = t^{nk} \cdot y$. The induced endomorphisms of
$A^{\fd}(*, M \times L)$ are also denoted by F$_n$ and are called
the \emph{Frobenius operations}.
\end{De}

To obtain the operations on the nil-terms we need to pass to the
operations on $\mathbb{D}_{\fd}(M \times L)$.
\begin{De}[Verschiebung on $\mathbb{D}_{\fd}(M \times L)$]
Define an exact functor
\begin{align*}
\text{V}_n \colon \mathbb{D}_{\fd}(M \times L) & \rightarrow \mathbb{D}_{\fd}(M \times L) \\
(Y_-,Y,Y_+) & \mapsto (Y_- \ot_{\varphi_n} \nn_-, Y \ot_{\varphi_n}
\zz, Y_+ \ot_{\varphi_n} \nn).
\end{align*}
\end{De}

\begin{De}[Frobenius on  $\mathbb{D}_{\fd}(M \times L)$]
Define an exact functor
\begin{align*}
\text{F}_n \colon \mathbb{D}_{\fd}(M \times L) & \rightarrow \mathbb{D}_{\fd}(M \times L) \\
(Y_-,Y,Y_+) & \mapsto (Y_{- \, \varphi_n},Y_{\varphi_n},Y_{+ \,
\varphi_n}).
\end{align*}
\end{De}
Recall the following \cite[Lemma 4.5]{fund1}
\begin{Le}\label{cong}
The forgetful functors
\begin{align*}
 \mathbb{D}_{\fd}(M \times \nn_+) & \rightarrow \mathbb{C}_{\fd}(M \times \nn_+) \\
 \mathbb{D}_{\fd}(M \times \nn_-) & \rightarrow \mathbb{C}_{\fd}(M \times \nn_-)
\end{align*}
induce homotopy equivalences $|h \mathcal{S}_\bullet
\mathbb{D}_{\fd}(M \times \nn_{\pm})| \rightarrow |h
\mathcal{S}_\bullet \mathbb{C}_{\fd}(M \times \nn_{\pm})|$.
\end{Le}
It is obvious that the map of Lemma \ref{cong} identifies the
operations defined on $\mathbb{D}_{\fd}(M \times L)$ with the
operations given on $\mathbb{C}_{\fd}(M \times L)$.
\begin{Pro}\label{restrict}
The Frobenius and Verschiebung operations restrict to operations on
$\NA(*,M)$ and $\NAM(*,M)$.
\end{Pro}
\begin{proof}
We start with the Verschiebung operations. Recall that the
augmentation map $\epsilon \colon \Omega |h \mathcal{S}_\bullet
\mathbb{D}_{\fd}(M \times \nn_+)| \rightarrow \Omega |h
\mathcal{S}_\bullet \mathbb{C}_{\fd}(M)|$ is induced by the functor
which maps an object $(Y_-,Y,Y_+)$ to $Y / \zz$. This implies that
the following diagram commutes:
\[
\xym{ \Omega |h \mathcal{S}_\bullet \mathbb{D}_{\fd}(M \times
\nn_+)| \ar[r]^-{\text{V}_n} \ar[d]_-{\epsilon} &
\Omega |h \mathcal{S}_\bullet \mathbb{D}_{\fd}(M \times \nn_+)| \ar[d]^-{\epsilon} \\
\Omega |h \mathcal{S}_\bullet \mathbb{C}_{\fd}(M)|
\ar[r]^-{\text{id}} & \Omega |h \mathcal{S}_\bullet
\mathbb{C}_{\fd}(M)|. }
\]
Since $\NA(*,M)$ is the homotopy fiber of $\epsilon$ we obtain that
the Verschiebung operation restricts to an operation on $\NA(*,M)$.

We continue with the Frobenius operations. Define an exact functor
I$ \colon \mathbb{C}_{\fd}(M) \rightarrow \mathbb{D}_{\fd}(M \times
\nn_+)$ by mapping $Y$ to $(Y \ot \nn_-, Y \ot \zz, Y \ot \nn_+)$
where $Y \ot L$ is defined to be $Y \times L / \sim$ where $(*,k)
\sim (*,l)$ for the base point $*$ and all $k,l$. Let $\iota$ denote
the map induced by I and notice that $\epsilon \circ \iota =
\text{id}$. This implies that $\NA(*,M)$ is the homotopy cofiber of
$\iota$. We also define a functor $\vee^n \colon \mathbb{C}_{\fd}(M)
\rightarrow \mathbb{C}_{\fd}(M)$ which maps $Y$ to $\vee^n Y$. The
Frobenius operation restricts to an operation on $\NA(*,M)$ since
the following diagram commutes:
\[
\xym{ \Omega |h \mathcal{S}_\bullet \mathbb{D}_{\fd}(M \times \nn_+)
| \ar[r]^{\text{F}_n}  &
\Omega |h \mathcal{S}_\bullet \mathbb{D}_{\fd}(M \times \nn_+) |  \\
\Omega |h \mathcal{S}_\bullet \mathbb{C}_{\fd}(M) | \ar[u]^{\iota}
\ar[r]^-{\vee^n} & \Omega |h \mathcal{S}_\bullet
\mathbb{C}_{\fd}(M)|. \ar[u]_-{\iota}}
\]
The proof that the operations restrict also to $\NAM(*,M)$ is
analogous.
\end{proof}
\subsection{Operations on the NIL-Category}
In this section we define Verschiebung and Frobenius on the
nil-category. This definition has the advantage that it is in
general easier to prove identities satisfied by the operations in
this setting.

\begin{De}[Verschiebung]
Define an exact functor
\begin{align*}
\text{V}'_n \colon \NIL_{\fd}(*,M) & \rightarrow \NIL_{\fd}(*,M)\\
(Y,f) & \mapsto
(\vee^n Y, \left( \begin{array}{cccc}
0         &        &           & f \\
\text{id} & \ddots &           &   \\
          & \ddots & 0         &   \\
          &        & \text{id} & 0
\end{array} \right)),
\end{align*}
where the matrix describes the map which sends the $n$-tuple $(y_0,
\ldots, y_{n-1})$ to $(f(y_{n-1}), y_0, \ldots , y_{n-2})$. The
induced endomorphisms of $K (\NIL_{\fd} (\ast,M))$ and on $\Nil
(\ast,M)$ are also denoted by $\text{V}'_n$.
\end{De}%
\begin{De}[Frobenius]
Define an exact functor
\begin{align*}
\text{F}'_n \colon \NIL_{\fd}(*,M) & \rightarrow \NIL_{\fd}(*,M) \\
(Y,f) & \mapsto (Y,f^{\circ n}).
\end{align*}
The induced endomorphisms of $K(\NIL_{\fd} (\ast,M))$ and on $\Nil
(\ast,M)$ are also denoted by $\text{F}'_n$.
\end{De}
%

%


\subsection{Operations on $K(\mathbb{S}[M \times L]$)}


The definition of the Verschiebung operations on the $\SS$-algebra
definition of the algebraic $K$-theory of spaces is needed for the
identification of the $\ZZ[\NN_\times]$-module structure on the
homotopy groups of the nil-terms in section \ref{sec:calculation}.
The definition of the Frobenius operations in this setup is not used
in the rest of the paper, but we include it for completeness.

Recall that we have the monoid $\SS$-algebra
\begin{equation*}
\SS[M \times L] \colon k_+ \mapsto (M \times L)_+ \wedge k_+ = M_+
\wedge L_+ \wedge k_+,
\end{equation*}
with the multiplication map defined via the monoid operation on $M
\times L$.

\begin{De}[Verschiebung]
Define an $\SS$-algebra map
\begin{align*}
V^{''}_n \colon \SS [M \times L] & \ra \SS[M \times L], \\
(m,l,i) & \mapsto (m,\varphi_n(l),i),
\end{align*}
where $m \in M$, $l \in L$ and $i \in k_+$. The induced map on
$K(\SS [M \times L])$ is also denoted by $V^{''}_n$ and is called
the {\it Verschiebung} operation.
\end{De}

The Frobenius operation is the corresponding transfer. Informally,
the map $\varphi_n$ makes the $\SS$-algebra $\SS[M \times L]$ into a
free $\SS[M \times L]$-module on $n$ generators. Using this any $m
\times m$ matrix over $\SS[M \times L]$ defines also an $mn \times
mn$ matrix over $\SS[M \times L]$. This assignment induces the
Frobenius operation. Recall that an $\AA$-module $\BB$ is a
$\Gamma$-space with an action map $\alpha \colon \AA \wedge \BB \ra
\BB$ satisfying certain relations, cf
\cite[1.5.1]{Dundas-Goodwillie-McCarthy(2004)}. Since the role of
the monoid $M$ in the sequel is not very interesting we omit it from
the notation and just talk about the $\SS$-algebra $\SS[L]$.

We denote by $\SS[L]_{\varphi_n}$ the $\SS[L]$ with the new
$\SS[L]$-module structure given by $\varphi_n$, i.e. with the action
map
\begin{equation*}
\begin{split}
\SS[L] \wedge \SS[L]_{\varphi_n} & \ra \SS[L]_{\varphi_n} \\
L_+ \wedge k_+ \wedge L_+ \wedge l_+ & \ra L_+ \wedge kl_+ \\
(z,i,u,j) & \mapsto (zn+u,ij).
\end{split}
\end{equation*}

With this $\SS[L]$-module structure we have an isomorphism of
$\SS[L]$-modules
\begin{equation*}
\begin{split}
n_+ \wedge \SS[L] & \ra \SS[L]_{\varphi_n} \\
n_+ \wedge L_+ \wedge k_+ & \ra L_+ \wedge k_+ \\
(j,z,i) & \ra (j+zn,i),
\end{split}
\end{equation*}
where we call $n_+ \wedge \SS[L]$ a free $\SS[L]$-module on $n$
generators.

Let $\AA$ be an $\SS$-algebra and let $\BB$, $\CC$ be $\AA$-modules.
There is a $\Gamma$-space of $\AA$-module maps from $\BB$ to $\CC$
defined as in \cite[Definition
1.5.4]{Dundas-Goodwillie-McCarthy(2004)}. Its $k$-th space can be
thought of as the space of natural transformations
\[
\lambda \in \Gamma \sS_\ast (\BB,\CC) (k_+) =
\Nat_{\Gamma\sS_\ast}(\BB(-),\CC(k_+ \wedge -)).
\]
such that the following diagram commutes
\[
\xymatrix{
\AA \wedge \BB (-) \ar[r] \ar[d]_{\id \wedge \lambda} & \BB (-) \ar[d]^{\lambda} \\
\AA \wedge \CC (k_+ \wedge -) \ar[r] & \CC (k_+ \wedge -). }
\]

In analogy with linear algebra there is an isomorphism $\AA \cong
\Map_\SS (\SS,\AA)$ of $\SS$-algebras, since an $\SS$-algebra map
$\SS \ra \AA$ is determined by what it does on $\SS(1_+) = 1_+$. It
follows that there are $\SS$-algebra isomorphisms
\[
\Mat_m (\AA) \; \cong \; \Map_{\SS} (m_+ \wedge \SS , m_+ \wedge
\AA) \; \cong \; \Map_{\AA} (m_+ \wedge \AA , m_+ \wedge \AA).
\]
The second isomorphism is given using the isomorphism $\SS \wedge
\AA \ra \AA$ and the multiplication $\AA \wedge \AA \ra \AA$ as the
composition
\[
\Map_\SS (m_+ \wedge \SS,m_+ \wedge \AA) \ra \Map_\AA (m_+ \wedge
\SS \wedge \AA,m_+ \wedge \AA \wedge \AA) \ra \Map_\AA (m_+ \wedge
\AA,m_+ \wedge \AA).
\]

A natural transformation $\SS[L](-) \ra \SS[L](k_+ \wedge -)$ which
represents a point in $\Map_{\SS[L]} (\SS[L],\SS[L])(k_+)$ is
determined by what it does on $-=1_+$ where it is of the form
\[
z \mapsto (z+a,i)
\]
for some $a \in L$ and $i \in k_+$. It follows that it also
represents a point in the space $\Map_{\SS[L]} (\SS[L]_{\varphi_n},
\SS[L]_{\varphi_n})(k_+)$. More generally we have an $\SS$-algebra
map
\[
\varphi_n^\ast \colon \Map_{\SS[L]} (m_+ \wedge \SS[L],m_+ \wedge
\SS[L]) \ra \Map_{\SS[L]} (m_+ \wedge \SS[L]_{\varphi_n},m_+ \wedge
\SS[L]_{\varphi_n}).
\]
Further we obtain an $\SS$-algebra map
\begin{multline*}
\varphi^\ast_n \colon \Mat_m (\SS[L]) \ra \Map_{\SS[L]} (m_+ \wedge
\SS[L],m_+ \wedge \SS[L]) \xrightarrow{\varphi_n^\ast} \\ \ra
\Map_{\SS[L]}
(m_+ \wedge \SS[L]_{\varphi_n},m_+ \wedge \SS[L]_{\varphi_n}) \ra \\
\ra \Map_{\SS[L]} (mn_+ \wedge \SS[L],mn_+ \wedge \SS[L]) \ra
\Mat_{mn} (\SS[L]).
\end{multline*}
This induces the maps of simplicial monoids $\widehat{\Mat}_m
(\SS[L]) \ra \widehat{\Mat}_{mn} (\SS[L])$ and
\begin{equation} \label{F_n}
\varphi^\ast_n \colon \widehat{\GL}_m (\SS[L]) \ra
\widehat{\GL}_{mn} (\SS[L]).
\end{equation}

\begin{De}[Frobenius]
The map induced by the maps $\varphi^\ast_n$ from (\ref{F_n}) is
denoted  $F^{''}_n \colon K (\SS[L]) \ra K(\SS[L])$  and called the
{\it Frobenius} operation.
\end{De}

\noindent The formula for the map $\varphi^\ast_n \colon \Mat_m
(\SS[L]) \ra \Mat_{mn} (\SS[L])$ in case $m = 1$ is
\[
\left( a \right) \quad \longmapsto \quad \left(\begin{array}{cccccc}
   &   &   & b \\
   &   &   &   & \ddots \\
   &   &   &   &   & b \\
b \\
   & \ddots \\
   &   & b
\end{array} \right)
\]
where $a = bn+c$, $0 \leq c < n$ and the first non-zero entry in the
first column is in the row $c+1$. For $m > 1$ replace any non-zero
entry $a$ by the matrix as above.

\begin{Rem}
In view of the definition of $K(\AA)$ we have defined the operations
$V^{''}_n$, $F^{''}_n$ only on the base point component. The
extension to the whole $K(\AA)$ can be done by artificially on $K_0
(\ZZ[\pi_0 \AA])$.
\end{Rem}


\section{Identifying Verschiebung and Frobenius}\label{indent}
In the previous section we defined the Frobenius and Verschiebung
operations using various definitions of the nil-terms. In this
section we prove that these definitions in fact coincide. Recall
from Theorem \ref{thm-fund3} that for a topological monoid $M$ there
are natural homotopy equivalences $\delta_\pm \colon \Nil (\ast,M)
\ra \Omega \NAPM (\ast,M)$.
\begin{Th}\label{identyfing}
The following diagrams commute up to a preferred homotopy
\[
\xymatrix{
\Nil (\ast,M) \ar[r]^{V'_n} \ar[d]_{\delta_\pm} & \Nil (\ast,M) \ar[d]^{\delta_\pm} \\
\Omega \NAPM (\ast,M) \ar[r]_{V_n} & \Omega \NAPM (\ast,M) } \quad
\xymatrix{
\Nil (\ast,M) \ar[r]^{F'_n} \ar[d]_{\delta_\pm} & \Nil (\ast,M) \ar[d]^{\delta_\pm} \\
\Omega \NAPM (\ast,M) \ar[r]_{F_n} & \Omega \NAPM (\ast,M), }
\]
for all $n \in \NN_\times$.
\end{Th}
For the proof we need the notion of the projective line category
\cite[page 4]{fund3}.
\begin{De}[Projective line category]
The \emph{projective line category} $\mathbb{P}_{\fd}(M)$ is given
by the full subcategory of $\mathbb{D}_{\fd}(M \times \zz)$ whose
objects $(Y_-,Y,Y_+)$ satisfy the additional property that the
induced maps $Y_-(t) \rightarrow Y$ and $Y_+(t) \rightarrow Y$ are
weak homotopy equivalences. A \emph{cofibration} in
$\mathbb{P}_{\fd}(M)$ is a morphism which is a cofibration in
$\mathbb{D}_{\fd}(M \times \zz)$. A \emph{weak homotopy equivalence}
from $(Y_-,Y,Y_+)$ to $(X_-,X,X_+)$ is a morphism in which $Y_-
\rightarrow X_-$, $Y \rightarrow X$ and $Y_+ \rightarrow X_+$ are
weak homotopy equivalences of spaces. Let
$\mathbb{P}^{h_{\nn_+}}_{\fd}(M)$ be the full Waldhausen subcategory
of $\mathbb{P}_{\fd}(M)$ with objects $(Y_-,Y,Y_+)$ such that $Y_+$
is acyclic.
\end{De}
Note that the operations of section \ref{sec:op-on-na} restrict to
$\mathbb{P}^{h_{\nn_+}}_{\fd}(M)$. Further we need a pair of exact
functors from \cite{fund3} which relate the projective line category
to the nil-category.
There is an exact functor
\begin{align*}
\Phi \colon \NIL_{\fd}(*,M) & \rightarrow
\mathbb{P}^{h_{\nn_+}}_{\fd}(M) \\
(Y,f) & \mapsto (Y_f,Y_f(t),*)
\end{align*}
where $Y_f$ is the homotopy coequalizer of the pair of maps
\[
\xym{Y \ot \nn_- \ar@<0.5ex>[r]^{f} \ar@<-0.5ex>[r]_{t^{-1}} & Y \ot
\nn_-}.
\]
where $t^{-1}$ is the map $(y,r) \mapsto (y,r-1)$. An explicit model
for the homotopy coequalizer of two maps $\alpha,\beta \colon U \ra
V$ is the quotient of $V \coprod (U \times [0,1])$ by $(u,0) \sim
\al(u)$, $(u,1) \sim \beta (u)$ and $* \times [0,1] \sim$  the base
point of $V$. Note that if $f$ is a cofibration then $Y_f$ is
homotopy equivalent to $Y$ with the $\NN_-$-action given by $(t^k,y)
\mapsto f^{\circ k}(y)$.

There is also an exact functor
\begin{align*}
\Psi \colon \mathbb{P}^{h_{\nn_+}}_{\fd}(M) & \rightarrow
\NIL_{\fd}(*,M) \\
(Y_-,Y,Y_+) & \mapsto (Y_-,t^{-1}).
\end{align*}

\begin{Le}\label{commu}
The following diagrams commute up to a preferred homotopy
\[
\xym{ \Omega |h \mathcal{S}_\bullet \NIL_{\fd}(*,M)| \ar[r]^-{\text{V}'_n} \ar[d]_-\Phi & \Omega |h \mathcal{S}_\bullet \NIL_{\fd}(*,M)| \ar[d]^-{\Phi} \\
  \Omega |h \mathcal{S}_\bullet \mathbb{P}^{h_{\nn_+}}_{\fd}(M)| \ar[r]^{\text{V}_n} & \Omega |h \mathcal{S}_\bullet \mathbb{P}^{h_{\nn_+}}_{\fd}(M)| }
\]
\[
\xym{\Omega |h \mathcal{S}_\bullet \NIL_{\fd}(*,M)| \ar[r]^-{\text{F}'_n} \ar[d]_-\Phi & \Omega |h \mathcal{S}_\bullet \NIL_{\fd}(*,M)| \ar[d]^-{\Phi} \\
     \Omega |h \mathcal{S}_\bullet \mathbb{P}^{h_{\nn_+}}_{\fd}(M)| \ar[r]^{\text{F}_n} & \Omega |h \mathcal{S}_\bullet \mathbb{P}^{h_{\nn_+}}_{\fd}(M)|.}
\]
for all $n \in \NN_\times$.
\end{Le}
\begin{proof}
We will first reduce to the case where $f$ is a cofibration.
Consider the functor $\Psi \circ \Phi \colon \NIL_{\fd}(*,M)
\rightarrow \NIL_{\fd}(*,M)$ which maps $(Y,f)$ to $(Y_f,t)$. By
construction we have that $t$ is a cofibration. Since $\Psi \circ
\Phi$ induces a homotopy equivalence \cite[Lemma 3.5]{fund3} we can
precompose our diagrams with $\Psi \circ \Phi$ and therefore assume
that all objects in $\NIL_{\fd}(*,M)$ are of the form $(Y,f)$ where
$f$ is a cofibration. Thus we can use for $Y_f$ the space $Y$ with
the $\nn_-$-action given via $f$.

To prove the lemma we need to define a preferred natural homotopy
between V$_n \circ \Phi$ and $\Phi \circ \text{V}_n$ and between
F$_n \circ \Phi$ and $\Phi \circ \text{F}_n$.

We start with the Verschiebung operations. We have
\[
\text{V}_n \circ \Phi \big((Y,f)\big)= (Y \ot_{\varphi_n} \NN_-,
Y(t) \ot_{\varphi_n} \zz,*)
\]
where $t^{-1}$ acts on $Y$ via $f$ and
\[
\Phi \circ \text{V}_n \big( (Y,f) \big) = ((\vee^n Y), (\vee^n Y)
(t),*)
\]
where $t^{-1}$ acts on $(\vee^n Y)$ via
\[
\left( \begin{array}{cccc}
0         &        &           & f \\
\text{id} & \ddots &           &   \\
          & \ddots & 0         &   \\
          &        & \text{id} & 0
\end{array} \right).
\]

Points in $Y(t) \ot_{\varphi_n} \zz$ are triples $(y,z,i)$ where $y
\in Y$, $0 \leq i \leq n-1$ and $z \in \zz$ modulo the equivalence
relation
\begin{align*}
(y,z,i) & \sim (f(y),z + 1,i)\\
(*,z,i) & \sim (*,z',i').
\end{align*}
The $\ZZ$-action is given by $t \cdot (y,z,i) \mapsto (y,z,i+1)$ for
$i \leq n-2$ and $t \cdot (y,z,n-1) \mapsto (y,z-1,1)$.

Points in $(\vee^n Y)(t)$ are triples $(y,i,z)$ where $y \in Y$, $0
\leq i \leq n-1$ and $z \in \zz$ modulo the equivalence relation
\begin{align*}
(y,i,z) & \sim (y,i+1,z+1) \,\,\, \text{for $i \neq n-1$} \\
(y,n-1,z) & \sim (f(y),0,z +1) \\
(*,i,z) & \sim (*,i',z').
\end{align*}
The $\ZZ$-action is given by $t \cdot (y,i,z) \mapsto (y,i,z - 1)$.

The map
\begin{align*}
Y(t) \ot_{\varphi_n} \zz & \rightarrow  (\vee^nY)(t) \\
(y,z,i) & \mapsto (y,0,z\cdot n + i)
\end{align*}
and a similar map $Y \ot_{\varphi_n} \NN_- \rightarrow  (\vee^nY)$
induce a natural equivalence of the functors and hence a preferred
homotopy between V$_n \circ \Phi$ and $\Phi \circ \text{V}_n$.

Now the Frobenius operations. We have
\[
\text{F}_n \circ \Phi \big((Y,f)\big)= (Y_{\varphi_n}
,Y(t)_{\varphi_n},*)
\]
where $t^{-1}$ acts on $Y$ via $f$ and
\[
\Phi \circ \text{F}_n \big( (Y,f) \big) = (Y ,Y (t),*)
\]
where now $t^{-1}$ acts on $Y$ via $f^{\circ n}$. But this is just
the space $Y_{\varphi_n}$, so we can choose the identity on the
first coordinate. It remains to find a preferred natural isomorphism
$Y_{\varphi_n} (t) \ra Y(t)_{\varphi_n}$.

Points in $Y(t)_{\varphi_n}$ are pairs $(y,z)$ where $y \in Y$ and
$z \in \zz$ modulo the equivalence relation
\[
(y,z) \sim (f(y),z - 1).
\]
The $\ZZ$-action is given by $t \cdot(y,z) \mapsto (y,z+n)$.

Points in $Y_{\varphi_n} (t)$ are pairs $(y,z)$ where $y \in Y$ and
$z \in \zz$ modulo the equivalence relation
\[
(y,z) \sim (f^{\circ n}(y),z - 1).
\]
The $\ZZ$-action is given by $t \cdot (y,z) \mapsto (y,z+1)$.

The map
\begin{align*}
 Y_{\varphi_n} (t) & \rightarrow Y(t)_{\varphi_n}\\
(y,z) & \mapsto (y,n\cdot z)
\end{align*}
induces a preferred natural equivalence of functors and hence a
preferred natural homotopy between F$_n \circ \Phi$ and $\Phi \circ
\text{F}_n$.
\end{proof}
\begin{proof}[Proof of Theorem \ref{identyfing}]
Let us briefly recall how are the terms $\NA$ and $\Nil$ identified
in \cite{fund3} (analogous proof works for $\NAM$). First it is
proven that the inclusion functors induce a homotopy fibration
sequence \cite[Proposition 2.2]{fund3}:
\[
\Omega | h \mathcal{S} _\bullet \mathbb{P}^{h_{\nn_+}}_{\fd}(M)|
\rightarrow \Omega |h \mathcal{S}_\bullet \mathbb{P}_{\fd}(M)|
\rightarrow  \Omega | h \mathcal{S}_\bullet \mathbb{D}_{\fd}(M
\times \nn _+ )|.\] Recall that $\NA$ is defined to be the homotopy
fiber of the augmentation map
\[
\epsilon \colon \Omega | h \mathcal{S}_\bullet \mathbb{D}_{\fd}(M
\times \nn _+ )| \rightarrow \Omega |h \mathcal{S}_\bullet
\mathbb{C}_{\fd}(M) |.
\]
We define $\Omega |h \mathcal{S}_\bullet \mathbb{P}_{\fd}(M)
|^\epsilon$ to be the homotopy fiber of the augmentation map
\[
\epsilon \colon \Omega |h \mathcal{S}_\bullet \mathbb{P}_{\fd}(M)|
\rightarrow \Omega |h \mathcal{S}_\bullet \mathbb{C}_{\fd}(M) |
\]
which is induced by the functor which maps $(Y_-,Y,Y_+)$ to $Y/
\zz$. Since
\[
\Omega |h \mathcal{S}_\bullet \NIL_{\fd}(*,M)| \simeq \Omega | h
\mathcal{S} _\bullet \mathbb{P}^{h_{\nn_+}}_{\fd}(M)|
\]
we obtain the following diagram where the vertical and horizontal
sequences of maps are homotopy fibration sequences
\[
\xym{& \Omega |h \mathcal{S}_\bullet \mathbb{P}_{\fd}(M) |^\epsilon
\ar@{-->}[r] \ar[d] &
\NA(*,M) \ar[d] \\
\Omega | h \mathcal{S}_\bullet \NIL_{\fd}(M,*) | \ar[r]
\ar@{-->}[ur]& \Omega | h \mathcal{S}_\bullet \mathbb{P}_{\fd}(M) |
\ar[r] \ar[d] &
\Omega | h \mathcal{S}_\bullet \mathbb{D}_{\fd}(M \times \nn _+ )| \ar[d] \\
& \Omega | h \mathcal{S}_\bullet \mathbb{C}_{\fd}(M) | \ar[r] &
\Omega | h \mathcal{S}_\bullet \mathbb{C}_{\fd}(M) |. }
\]
This implies that we obtain the dashed homotopy fibration sequence. The map
\[
\Omega |h \mathcal{S}_\bullet \mathbb{P}_{\fd}(M) |^\epsilon
\rightarrow \NA(*,M)
\]
is null homotopic and there is a homotopy equivalence
\[
\Omega |h \mathcal{S}_\bullet \mathbb{P}_{\fd}(M) |^\epsilon \simeq
\Omega |h \mathcal{S}_\bullet \mathbb{C}_{\fd}(M) |.
\]
All these constructions are natural in $M$ and hence we obtain a
natural homotopy equivalence $\delta_+ \colon \Nil(*,M)
\xrightarrow{\simeq} \Omega \NA(*,M)$.

Define the Verschiebung operations on $\Omega |h \mathcal{S}_\bullet
\mathbb{C}_{\fd}(M)|$ to be the identity and define the Frobenius
operation on $\Omega |h \mathcal{S}_\bullet \mathbb{C}_{\fd}(M)|$ to
be the map which is induced by the functor $\vee^n$ defined in the
proof of Proposition \ref{restrict}.

We prove the result only for the Verschiebung operations. The proof
for the Frobenius operation is analogous with $\epsilon$ replaced by
$\iota$. We have the following diagram

{\footnotesize
\begin{center}
 \makebox[0pt]{
\xym@C=0.1em@R=1em{& & & \Omega |h \mathcal{S}_\bullet
\mathbb{P}_{\fd}(M) |^\epsilon \ar[dd] & &
\NA(*,M) \ar[dd] \\
&& \Omega |h \mathcal{S}_\bullet \mathbb{P}_{\fd}(M) |^\epsilon
\ar[ru]^{V_n} \ar[dd] &&
\NA(*,M) \ar[dd] \ar[ru]^{V_n} & \\
& \Omega | h \mathcal{S}_\bullet \NIL_{\fd}(*;M) |
\ar[rr]|(.53){\phantom{**}}  && \Omega | h \mathcal{S}_\bullet
\mathbb{P}_{\fd}(M) |  \ar[rr]|(.465){\phantom{**}}
\ar[dd]|(.49){\phantom{M}}  &&
\Omega | h \mathcal{S}_\bullet \mathbb{D}_{\fd}(M \times \nn _+ )| \ar[dd] \\
\Omega | h \mathcal{S}_\bullet \NIL_{\fd}(*;M) | \ar[rr]
\ar[ru]^{V_n'} && \Omega | h \mathcal{S}_\bullet \mathbb{P}_{\fd}(M)
| \ar[ru]^{V_n} \ar[rr] \ar[dd] &&
\Omega | h \mathcal{S}_\bullet \mathbb{D}_{\fd}(M \times \nn _+ )| \ar[ru]^{V_n} \ar[dd] & \\
& & & \Omega |h \mathcal{S}_\bullet \mathbb{C}_{\fd}(M) |
\ar[rr]|(.47){\phantom{**}} &&
\Omega |h \mathcal{S}_\bullet \mathbb{C}_{\fd}(M) | \\
& & \Omega |h \mathcal{S}_\bullet \mathbb{C}_{\fd}(M) |
\ar[ru]^{V_n} \ar[rr] && \Omega |h \mathcal{S}_\bullet
\mathbb{C}_{\fd}(M) | \ar[ru]^{V_n} &. }}
\end{center}
}

The vertices in the top layer of the diagram are defined as homotopy
fibers of the corresponding vertical maps. The horizontal sequences
are homotopy fibration sequences. After using the trick of
precomposing with $\Psi \circ \Phi$ (as in Lemma \ref{commu}) we
obtain by the proof of Proposition \ref{restrict} that the diagram
is strictly commutative except at the left horizontal square where
it is commutative only up to a preferred homotopy by Lemma
\ref{commu}. This diagram can be replaced by a (preferred) homotopy
equivalent diagram which is strictly commutative. Adding the dashed
arrows explained above we obtain a self-map of the dashed homotopy
fibration sequence (consisting of three maps so that the
corresponding squares strictly commute). The statement of theorem
then follows since the natural homotopy equivalence $\delta_+$ is
obtained from the ``connecting'' map of the dashed homotopy
fibration sequence.
\end{proof}
If $M$ is a topological monoid, there exists a natural homotopy
equivalence $\varepsilon \colon A^{\fd} (\ast, M)
\xrightarrow{\simeq} K(\SS[M])$. This was mentioned in Remark
\ref{ident-A-K} when $M = \Omega X$ and for general $M$ the
construction of $\varepsilon$ is sketched in the proof of the
following theorem. 

\begin{Th}
Let $L = \NN_+$, $\NN_-$ or $\ZZ$. Then the following diagrams commute up to a preferred homotopy
\[
\xymatrix{
A^{\fd} (\ast,M \times L) \ar[r]^{V_n} \ar[d]_{\varepsilon} & A^{\fd} (\ast,M \times L) \ar[d]^{\varepsilon} \\
K(\SS[M \times L]) \ar[r]_{V_n''} & K(\SS[M \times L]) } \quad
\xymatrix{
A^{\fd} (\ast,M \times L) \ar[r]^{F_n} \ar[d]_{\varepsilon} & A^{\fd} (\ast,M \times L) \ar[d]^{\varepsilon} \\
K(\SS[M \times L]) \ar[r]_{F_n''} & K(\SS[M \times L]) }
\]
for all $n \in \NN_\times$.
\end{Th}

\begin{proof}
Recall briefly the definition of $\varepsilon \colon A^{\fd} (\ast,M
\times L) \ra K(\SS[M \times L])$ using the notation of
\cite{Waldhausen1126} and \cite{Dundas-Goodwillie-McCarthy(2004)}.
There $A^{\fd} (\ast,M \times L)$ is defined by the
$\sS_\bullet$-construction on the category $R_{\fd}(*,M \times L)$,
of all finitely-dominated objects which has the structure of a
Waldhausen category in a similar way as $\mathbb{C}_{\fd}(*,M)$. The
required identification goes via a sequence of the following weak
equivalences natural in $M \times L$ \cite[page
385-389]{Waldhausen1126}, \cite[page
113]{Dundas-Goodwillie-McCarthy(2004)}:
\begin{align*} \Omega |h
\mathcal{S}_{\bullet} R_{fd}(*,M \times L)|
& \simeq \Omega |\colim_n h \mathcal{N}_{\bullet} R^n_{\fd}(*,M \times L)| \\
& \simeq K_0 (\ZZ[\pi_0 \AA]) \times \colim_{n,k} |h R^n_{k}(*,M \times L)|^+ \\
 & \simeq K_0 (\ZZ[\pi_0 \AA]) \times \colim_{n,k} B \mathcal{H}^n_k (M \times L)^+ \\
 & \simeq K_0 (\ZZ[\pi_0 \AA]) \times \colim_{k} B \widehat{\GL}_k(\mathbb{S}[M \times L])^+
\end{align*}
The first line is Proposition 2.2.2 in \cite{Waldhausen1126}
(slightly modified, see also Remark on page 389), the second line is
Segal's group completion theorem \cite{Segal}, the third line is
Proposition 2.2.5 in \cite{Waldhausen1126} and the last line is
explained in \cite[page 113]{Dundas-Goodwillie-McCarthy(2004)}.

Now the identification of the Verschiebung operations follows from
naturality of the above weak equivalences and the fact that in both
cases the Verschiebung is induced by the monoid morphism $\varphi_n
\colon M \times L \ra M \times L$.

The identification of the Frobenius operations is a tedious
verification. Informally one defines the Frobenius operations in all
the lines of the above identification, on the first two lines they
are defined similarly as on $\mathbb{C}_{\fd}(*,M \times L)$ and on
the third line they are defined similarly as on
$\widehat{\GL}_k(\mathbb{S}[M \times L])^+$. The details are left
for the reader.
\end{proof}
\section{Identities}
In this section we prove certain useful identities satisfied by the
Frobenius and Verschiebung operations. As a consequence we obtain
Corollary \ref{infinite} about the non-finiteness generation of the
nil-terms. Further recall that the nil-terms are infinite loop
spaces, so we can think of them as of spectra and hence as of
$\SS$-modules. Using this and some of the identities proven here we
obtain Corollaries \ref{ZZ[NN]-module}, \ref{ZZ[NN*NN]-module} about
the structure of an $\SS[\NN_\times]$-module and an $\SS[\NN_\times
\ast \NN_\times]$-module on the nil-terms and hence a structure of a
$\ZZ[\NN_\times]$-module and a $\ZZ[\NN_\times \ast
\NN_\times]$-module on the homotopy groups of the nil-terms. Here
$\NN_\times = \{1,2,\ldots\}$ is considered as a multiplicative
monoid and $\NN_\times \ast \NN_\times$ denotes the free product of
monoids.

\begin{Pro}\label{Relation1}
For $i \geq 0$ and $m,n \in \NN_\times$ the following identities
hold on $\NA (\ast,M)$:
\begin{align*}
\text{F}_1 = \text{V}_1 & = 1 \\
\text{V}_n \text{V}_m & = \text{V}_{n \cdot m} \\
\text{F}_n \text{F}_m & = \text{F}_{n \cdot m} \\
\text{F}_n \text{V}_n (x) & \simeq n \cdot x \\
\text{F}_n \text{V}_m & \simeq  \text{V}_m \text{F}_n \; \; \textrm{if} \; \; m,n \; \; \textrm{are coprime.}\\
\end{align*}
\end{Pro}
\begin{proof}
The first identity is obvious, the proof of the next two follows
from the formula $\varphi_{mn} = \varphi_n \circ \varphi_m$. For the
fourth identity let $(Y,f)$ be an object in $\NIL_{\fd}(*;M)$. Then
matrix manipulation yields
\begin{align*}
\text{F}_n \text{V}_n \big( (Y,f)\big)
&  = \text{F}_{n} \Big((\vee^n Y,\text{%
\fontsize{8pt}{6pt}\selectfont $ \left(
\begin{array}{cccc} 0 & & & f
\\ \text{id} & \ddots & &\\ & \ddots & 0 & \\ & & \text{id} & 0 \end{array}
 \right) $} \Big) \\
& \cong \Big(\vee^n Y,\text{%
\fontsize{8pt}{6pt}\selectfont $ \left(
\begin{array}{cccc}
f &        &   & \\
& f        &   & \\
& & \ddots &     \\
& &        & f
\end{array} \right)$} \Big) \cong \vee^n (Y,f),
\end{align*}
where the $\cong$-signs mean that the corresponding functors are
isomorphic. The resulting map induces multiplication by $n$ on
$\NIL_{\fd}(*;M)$ and therefore also on $\NA (\ast,M)$.

For the last identity we will show that there is a cofibration
sequence of functors
\begin{equation} \label{cofib-seq}
\xym{\text{V}_m  \text{F}_n  \Psi \circ \Phi \ar[r]^{\Xi} &
\text{F}_n \text{V}_m \Psi \circ \Phi \ar[r] & \text{G}},
\end{equation}
where $\Phi, \Psi$ are the functors from section \ref{indent} and G
is functor which induces the trivial map on $\Nil(*,M)$. The fact
that $\Psi \circ \Phi \simeq \text{id}$ and Waldhausen additivity
theorem \cite[Proposition 1.3.2(4)]{Waldhausen1126} now imply the
identity. Recall that the effect of $\Psi \circ \Phi$ on an object
$(Y,f)$ in $\NIL_{\fd} (\ast, M)$ is that it turns the map $f$ to a
cofibration, so we may assume that $f$ is a cofibration. To obtain
the cofibration sequence (\ref{cofib-seq}) note that both values
$\text{V}_m  \text{F}_n (Y,f)$ and $\text{F}_n \text{V}_m (Y,f)$ are
of the shape $(\vee_m Y,?) = (Y_1 \vee \ldots \vee Y_m, ?)$, where
the maps $?$ permute all $Y_i$. More precisely, since $m$ and $n$
are coprime the maps $?$ are of the following shape:
\[
\xym{ \text{V}_m  \text{F}_n && \text{F}_n \text{V}_m \\
      Y_1 \ar[d]_{\text{id}} && Y_{\sigma(1)} \ar[d]_{f^{l_{\sigma(1)}}} \\
      Y_2 \ar[d]_{\text{id}} && Y_{\sigma(2)} \ar[d]_{f^{l_{\sigma(2)}}} \\
      Y_3 \ar[d]_{\text{id}} && Y_{\sigma(3)} \ar[d]_{f^{l_{\sigma(3)}}} \\
      \vdots \ar[d]_{\text{id}} && \vdots \ar[d]_{f^{l_{\sigma(m-1)}}} \\
      Y_m \ar@/^3pc/[uuuu]^{f^n} && Y_{\sigma(m)}  \ar@/_3pc/[uuuu]_{f^{l_{\sigma(m)}}} }
\]
where $\sigma$ is a permutation on $m$ letters and $\sum_i
l_{\sigma(i)} = n$. Define the natural transformation:
\[
\xym{ \text{V}_m  \text{F}_n \ar[rr]^\Xi  && \text{F}_n \text{V}_m \\
      Y_1 \ar[d]_{\text{id}} \ar[rr]^{\text{id}} && Y_{\sigma(1)} \ar[d]_{f^{l_\sigma(1)}} \\
      Y_2 \ar[d]_{\text{id}} \ar[rr]^{f^{l_{\sigma(1)}}} && Y_{\sigma(2)} \ar[d]_{f^{l_{\sigma(2)}}} \\
      Y_3 \ar[d]_{\text{id}} \ar[rr]^{f^{l_{\sigma(1)} + l_{\sigma(2)}}}     && Y_{\sigma(3)} \ar[d]_{f^{l_{\sigma (3)}}} \\
      \vdots \ar[d]_{\text{id}} && \vdots \ar[d]_{f^{l_{\sigma(m-1)}}} \\
      Y_m \ar@/^3pc/[uuuu]^{f^n} \ar[rr]^{f^{n - l_{\sigma(m)}}} && Y_{\sigma(m)}
      \ar@/_3pc/[uuuu]_{f^{l_{\sigma(m)}}}}
\]
and the functor G as a cofiber of $\Xi$. An inspection shows that
the morphism of G$(Y,f)$ is a strictly lower triangular matrix and
therefore is trivial in $\Nil (*,M)$.
\end{proof}

\begin{Co} \label{ZZ[NN]-module}
The Verschiebung operations $V_n$ for $n \in \NN_\times$ define a
structure of an $\SS[\NN_\times]$-module on $\NAPM (\ast,M)$ and
hence a structure of a $\ZZ[\NN_\times]$-module on $\pi_i \NAPM
(\ast,M)$.
\end{Co}

The statement follows from the identity $\text{V}_n \text{V}_m =
\text{V}_{n \cdot m}$ of Proposition \ref{Relation1}.
\begin{Co} \label{ZZ[NN*NN]-module}
The Frobenius and Verschiebung operations $F_n$, $V_n$ for $n \in
\NN_\times$ define a structure of an $\SS[\NN_\times \ast
\NN_\times]$-module on $\NAPM (\ast,M)$ and hence a structure of a
$\ZZ[\NN_\times \ast \NN_\times]$-module on $\pi_i \NAPM (\ast,M)$.
\end{Co}
The statement follows from the identities $\text{F}_n \text{F}_m =
\text{F}_{n \cdot m}$, $\text{V}_n \text{V}_m = \text{V}_{n \cdot
m}$ of Proposition \ref{Relation1}.

Now we are ready for the proof of Theorem \ref{infinite}.

\begin{proof}[Proof of Theorem \ref{infinite}]
Since F$_m$ and V$_m$ satisfy the relation stated in
Theorem~\ref{Relation1} and $\text{F}_m (x) = 0$ for $m$ bigger than
a certain number $M$, we can apply a trick which is due to Farrell
\cite{Farrell}.

For $n \in \nn$ let $h \mathcal{S}^n_{\bullet} \NIL(*;M)$ be the
full subcategory of $h \mathcal{S}_{\bullet} \NIL_{\fd}(*;M)$
consisting of objects of nilpotency degree smaller or equal to $n$.
We have
\begin{align*}
K_i\big(\NIL_{\fd}(*;M) \big) & = \pi_i \Omega |h \mathcal{S}_{\bullet} \NIL_{\fd}(*;M)| \\
                           & = \colim_n \pi_i \Omega |h \mathcal{S}_{\bullet}^n \NIL_{\fd}(*;M)|
\end{align*}
and
\begin{align*}
K_i\big( \mathbb{C}_{\fd}(M) \big) = \pi_i \Omega |h
\mathcal{S}_{\bullet}^0 \NIL_{\fd}(*;M)|.
\end{align*}

Assume now that $\pi_i \NA (\ast,M)$ is a finitely generated abelian
group. Thus we can find an $N$ such that the generators of $K_i
\big( \NIL_{\fd}(*; M) \big)$ are contained in $\pi_i \Omega|h
\mathcal{S}_{\bullet}^N \NIL_{\fd}(*;M)|$. This implies that there
is an $\ell$ such that $\text{F}_m$ is the trivial map for $m \geq
\ell$. Let $T$ be the torsion subgroup of $\pi_i \NA (\ast,M)$.
Since $\pi_i \NA (\ast,M)$ is finitely generated we have that $|T|$
is finite. Choose $t \in \nn$ such that $t \cdot |T| + 1 \geq \ell$.
By Theorem \ref{Relation1} we get that $\text{F}_{t \cdot |T| + 1}
\text{V}_{t \cdot |T| + 1}$ is a monomorphism. On the other hand,
since $t \cdot |T| + 1 \geq \ell$, the group $\pi_i \NA (\ast,M)$ is
in the kernel of $\text{F}_{t \cdot |T| + 1}$. Thus $\pi_i \NA
(\ast,M)$ is the trivial group.

The proof for $p$-primary subgroups is identical.
\end{proof}
%


\section{A calculation via Trace Invariants}
\label{sec:calculation}


This section contains the proof of Theorem \ref{pi-Nil-calculation}.
It is a calculation of the homotopy type of the $p$-completion of
the spaces $\NAPM (\ast)$ and their homotopy groups as
$\ZZ_p[\NN_\times]$-modules, where $p$ is an odd prime and $\ZZ_p$
denotes $p$-adic integers, in a certain range. We follow a general
scheme of Madsen \cite[section 4.5]{Madsen(1995)} to study the
linearization map $A^{\fd}(S^1) \ra K(\ZZ[\ZZ])$. The idea of the
calculation is based on a calculation of the linearization map
$A(\ast) \ra K(\ZZ)$ by B\"okstedt and Madsen
\cite{Boekstedt-Madsen(1994)}, who in turn claim the original idea
to come from Goodwillie.

Throughout we fix an odd prime $p$ and we denote the $p$-completion
of a space or of a spectrum by $\pco$. For details on $p$-completion
of spectra see \cite[Appendix
A.1.11]{Dundas-Goodwillie-McCarthy(2004)}. We note that connective
spectra are $p$-good (the $p$-completion map induces an isomorphism
on $\widetilde{H}_\ast (-; \FF_p)$) and the $p$-completion commutes
with fibration and cofibration sequences of spectra. The homotopy
groups of a $p$-completion of a spectrum become $\ZZ_p$-modules and
the Verschiebung operations induce a structure of a
$\ZZ_p[\NN]$-module on $\pi_\ast \NAPM (X)\pco$.

The scheme for calculating the homotopy type of $\NA (\ast)$ as
suggested by Madsen at the end of \cite[section 4.5]{Madsen(1995)}
is as follows. The product $\NA(\ast) \times \NAM (\ast)$ is the
homotopy cofiber of the assembly map $S^1_+ \wedge A(\ast) \ra
A(S^1)$. The product $\NK_+ (\ZZ) \times \NK_- (\ZZ)$ is the
homotopy cofiber of the assembly map $S^1_+ \wedge K(\ZZ) \ra
K(\ZZ[\ZZ])$ and is known to be contractible. For a space $X$,
denote by $F_K(X)$ the homotopy fiber of the linearization map $l
\colon A(X) \ra K(\ZZ[\pi_1 X])$ and recall that the two assembly
maps are compatible with $l$ and so we obtain a homotopy fibration
sequence
\begin{equation}\label{lin-assembly}
S^1_+ \wedge F_K (\ast) \ra F_K (S^1) \ra \NA (\ast) \times \NAM
(\ast).
\end{equation}
So it is enough to study $F_K (S^1)$, factor out the image of the
assembly map and identify the two summands. This is the approach we
take, we describe the space $F_K (S^1)\pco$ in a certain range in
Proposition \ref{F-K-S^1-calculation}. The proof of Theorem
\ref{pi-Nil-calculation} assuming \ref{F-K-S^1-calculation} is given
at the end of this section as well as the proofs of its corollaries.
In fact our calculation applies more generally, so we determine the
homotopy type of the space $F_K (B\pi)\pco$ for a discrete group
$\pi$ in a certain range.

\begin{Pro}\label{F-K-Bpi-calc}
Let $p$ be an odd prime and let $\pi$ be a discrete group. Then
there is a $(4p-7)$-connected map
\[
 \Big(\bigvee_{[g] \in [\pi]} \Sigma^{2p-2} \HH \FF_p \wedge (B\pi_+) \Big)\pco \lra
F_K (BZ_{\pi}(g))\pco~,
\]
where $[\pi]$ denotes the set of conjugacy classes of elements of
$\pi$ and $Z_{\pi} (g)$ is the centralizer of an element $g \in \pi$
.
\end{Pro}

\begin{Co} \label{Wh(Bpi)}
Let $p$ be an odd prime and let $\pi$ be a discrete group such that
the assembly map $ B\pi_+ \wedge K(\ZZ) \ra K(\ZZ \pi) $ is a
homotopy equivalence. Then there is a $(4p-7)$-connected map
\[
 \Big(\bigvee_{[g] \in [\pi] \smallsetminus [1]} \Sigma^{2p-2} \HH \FF_p \wedge (BZ_{\pi}(g)_+) \Big)\pco \lra \Wh (B\pi)\pco
\]
where $\Wh (B \pi)$ denotes the homotopy cofiber of the assembly map
$ B\pi_+ \wedge A(\ast) \ra A(B\pi)$.
\end{Co}

If $\pi = \ZZ$ we obtain the result for $F_K (S^1)$ and $\Wh (S^1) =
\NA(\ast) \times \NAM (\ast)$, since $S^1 = B \ZZ$. In this case we
are also interested in the Frobenius and Verschiebung operations.
Recall that these were defined also on $A^{\fd} (S^1)$. They
restrict to $F_K (S^1)$ and so we have a structure of a $\ZZ_p
[\NN_\times]$-module and a $\ZZ_p[\NN_\times \ast
\NN_\times]$-module on the $p$-completion of $F_K (S^1)$. Below
$S^1(n)$ denotes just a copy of $S^1$ indexed by $n \in \ZZ$.

\begin{Pro} \label{F-K-S^1-calculation}
If $p$ is an odd prime the following holds:
\begin{enumerate}
\item There is a $(4p-7)$-connected map
\[
\bigvee_{n \in \ZZ} \Sigma^{2p-2} \HH \FF_p \wedge (S^1(n)_+) \lra
F_K (S^1)\pco~,
\]
hence
\begin{align*}
\pi_{2p-2} F_K(S^1)\pco & \cong \oplus_{n \in \ZZ} \FF_p \{\beta_{n}\} \\
\pi_{2p-1} F_K(S^1)\pco & \cong \oplus_{n \in \ZZ} \FF_p \{\gamma_{n}\} \\
\pi_i F_K(S^1)\pco & \cong 0 \; \textup{for } i < 2p-2 \textup{ , }
2p-1 < i \leq 4p-7
\end{align*}
where
\begin{align*}
\beta_{n} & \in \pi_{2p-2} \Sigma^{2p-2} \HH \FF_p \wedge
(S^1(n)_+) \cong \FF_p \\
\gamma_{n} & \in \pi_{2p-1} \Sigma^{2p-2} \HH \FF_p \wedge (S^1(
n)_+) \cong \FF_p
\end{align*}
represent a certain choice of generators of these $\ZZ_p$-modules.
\item The $\ZZ_p[\NN_\times]$-module structure on $\pi_\ast F_K(S^1)\pco$
is given by
\begin{align*}
(n , \beta_m) & \mapsto \beta_{nm} \\
(n , \gamma_m) & \mapsto n \cdot \gamma_{nm}.
\end{align*}
\end{enumerate}
\end{Pro}

\begin{Co} \label{ZZ[NN]-module-F_K}
The  $\ZZ_p[\NN_\times]$-module $\pi_{2p-2} F_K (S^1)\pco$ is
generated by the elements $\beta_{-1}, \beta_0, \beta_1$, and the
$\ZZ_p[\NN_\times]$-module $\pi_{2p-1} F_K (S^1)\pco$ is not
finitely generated.
\end{Co}

\begin{Co} \label{ZZ[NN*NN]-module-F_K}
The  $\ZZ_p[\NN_\times \ast \NN_\times]$-module $\pi_{2p-1} \NAPM
(\ast,M)\pco$ is not finitely generated.
\end{Co}

A standard strategy when studying $K(\AA)$ for an $\SS$-algebra
$\AA$ is to use invariants of $K$-theory such as topological
Hochschild homology $\THH (\AA)$ or topological cyclic homology $\TC
(\AA,p)$ for a prime $p$, which we think of as functors from
$\SS$-algebras to spectra. This approach is particularly convenient
in the relative situation, and general results can be found for
example in \cite{Dundas-Goodwillie-McCarthy(2004)},
\cite{Madsen(1995)}. However, calculations for an arbitrary
$\SS$-algebra $\AA$ are still very hard. In the special case, when
$\AA = \SS [\Omega X]$, there is a convenient relationship between
$K$-theory, the above invariants and the free loop space of $X$ (see
Theorem \ref{Goodwillie-Calc-I} below). Hence, a possible approach
to our problem is to replace the linearization map by a map between
$K$-theory of the $\SS$-algebras of the form $\SS[\Omega X]$, which
in turn can be understood via the free loop space. It turns out that
this works if we restrict ourselves to a certain dimension range.
This is the already mentioned trick used by B\"okstedt and Madsen
\cite{Boekstedt-Madsen(1994)}.

We prove the following proposition which replaces the space $F_K (B
\pi)\pco$ for a discrete group by another (more approachable) space.
For a space $X$ we denote $Q (X) := \Omega^\infty \Sigma^\infty X$
and we let $SG = Q(S^0)_1$ be the identity component which is a
topological monoid with respect to the composition product.

\begin{Pro} \label{K-thy-trick-BM}
If $\pi$ is a discrete group, $p$ an odd prime, then there is a
$(4p-7)$-cartesian square
\[
\xymatrix{
A(BSG \times B\pi)\pco \ar[r] \ar[d] & A(B\pi)\pco \ar[d] \\
A(B\pi)\pco \ar[r] & K (\ZZ [\pi])\pco, }
\]
where the left vertical arrow is induced by the projection map and
the right vertical arrow is the linearization map.
\end{Pro}

Hence in the range $\ast \leq 4p-7$ it is enough to study the
homotopy fiber of the left vertical map. The proof of Proposition
\ref{K-thy-trick-BM} is given later. We first prove Proposition
\ref{F-K-Bpi-calc} assuming \ref{K-thy-trick-BM}.

We denote $\Lambda X = \Map (S^1,X)$, $\ev \colon \Lambda X \ra X$
is the evaluation map $\ev (\lambda) = \lambda (1)$. For a map $f
\colon Y \ra X$ the symbol $\Lambda (f \colon Y \ra X)$ denotes the
pullback of the diagram
\[
\xymatrix{ Y \ar[r]^{f} & X & \Lambda X \ar[l]_{\ev}.  }
\]
Goodwillie constructs in \cite{Goodwillie(1990)} a natural map $A(X)
\ra \Sigma^\infty_+ \Lambda X$. Given a map of spaces $f \colon Y
\ra X$ the composition $A(Y) \ra A(X) \ra \Sigma^\infty_+ \Lambda X$
factors through $\Sigma^\infty_+ \Lambda (f \colon Y \ra X)$. This
factorization has the following property, see \cite[Corollary
3.3]{Goodwillie(1990)}.

\begin{Th}[Goodwillie]
\label{Goodwillie-Calc-I} If $f \colon Y \ra X$ is a $k$-connected
map, then the square
\[
\xymatrix{ A(Y) \ar[r] \ar[d] & \Sigma^\infty_+ \Lambda (f \colon Y
\ra X)
\ar[d] \\
A(X) \ar[r] & \Sigma^\infty_+ \Lambda X }
\]
is $(2k-1)$-cartesian.
\end{Th}

\begin{proof}[Proof of Proposition \ref{F-K-Bpi-calc}]
By Theorem \ref{K-thy-trick-BM} we have a $(4p-7)$-connected map
\[
\hofiber \big( A(BSG \times B \pi) \ra A(B \pi) \big)\pco \ra F_K (B
\pi)\pco
\]
whose source can be understood using Proposition
\ref{Goodwillie-Calc-I}. If $f \colon Z \times X \ra X$ is the
projection, then $\Lambda (f  \colon Z \times X \ra X) = Z \times
\Lambda X$. If $Z$ has a base point the vertical maps in the diagram
of Theorem \ref{Goodwillie-Calc-I} have sections and the homotopy
fibers of these vertical maps are homotopy equivalent to the
homotopy cofibers of the corresponding sections. Thus, if $Z$ is a
$k$-connected space, we obtain a $(2k-1)$-connected map
\begin{equation} \label{fiber1}
\hofiber (A(Z \times X) \ra A(X) ) \ra \Sigma^\infty Z \wedge
(\Lambda X_+).
\end{equation} (We have a homotopy cofibration
sequence $\ast \times \Lambda X \ra Z \times \Lambda X \ra Z \wedge
(\Lambda X_+)$.)

For a discrete group $\pi$ we have a decomposition
\[
\Lambda B\pi \simeq \coprod_{[g] \in [\pi]} B Z_{\pi}(g),
\]
where $[\pi]$ denotes the set of conjugacy classes of elements in
$\pi$ and $Z_{\pi} (g)$ is the centralizer of an element $g \in
\pi$. Further, by \cite{Toda(1958)} we have a $(4p-6)$-connected map
\[
BSG\pco \ra \Sigma^{2p-2} \HH \FF_p.
\]
Plugging $Z = BSG$, $X = B\pi$ in (\ref{fiber1}) and $p$-completing
yields the statement of Proposition \ref{F-K-Bpi-calc}.
\end{proof}

Now we prove Proposition \ref{K-thy-trick-BM}. As already mentioned
B\"okstedt and Madsen in \cite[Proposition
9.11]{Boekstedt-Madsen(1994)} show such a result for the functor
$\TC(-,p)$ for an odd prime $p$ instead of $K$-theory, and with
$\pi$ the trivial group. This is our starting point. We will not
recall the definition of $\TC(-,p)$ here, since it is considerably
complicated. We need to know that for a prime $p$ there is a functor
from $\SS$-algebras to spectra $\AA \mapsto \TC(\AA,p)$ which comes
with a natural transformation $\trc \colon K(\AA) \ra \TC (\AA,p)$,
called the {\it cyclotomic trace map} which has certain properties.
A theorem of Dundas says that under certain assumptions relative
$K$-theory equals relative $\TC$ (Theorem \ref{Dundas} below).
Further we need another property of $\TC$ which we also only cite.

\begin{Th}[Dundas] \label{Dundas}
If $\AA \ra \BB$ is a map of $\SS$-algebras inducing an isomorphism
$\pi_0 (\AA) \ra \pi_0 (\BB)$, then the following square
\[
\xymatrix{ K(\AA)\pco \ar[r] \ar[d] & K(\BB)\pco \ar[d] \\
\TC (\AA,p)\pco \ar[r] & \TC (\BB,p)\pco. }
\]
is homotopy cartesian.
\end{Th}

For the proof see \cite{Dundas(1997)}, \cite{Madsen(1995)}.

\begin{Pro} \label{TC-properties}
The functor $\AA \mapsto \TC (\AA,p)$ has the following property. If
\[
\xymatrix{
\AA \ar[d] \ar[r] & \BB \ar[d] \\
\CC \ar[r] & \DD }
\]
is a $k$-cartesian square of connective $\SS$-algebras, then
\[
\xymatrix{
\TC(\AA,p) \ar[d] \ar[r] & \TC(\BB,p) \ar[d] \\
\TC(\CC,p) \ar[r] & \TC(\DD,p)
 }
\]
is a $(k-1)$-cartesian square.
\end{Pro}

For the proof combine the statements \cite[Proposition
1.4.2]{Dundas-Goodwillie-McCarthy(2004)} and \cite[Proposition
2.2.7]{Dundas-Goodwillie-McCarthy(2004)}. We also need the following
result from homotopy theory.

\begin{Le} \label{square1}
There is a $(4p-6)$-cartesian square of $\SS$-algebras
\begin{equation} \label{square1-diagram}
\begin{split}
\xymatrix{ \SS[SG \times \pi]\pco \ar[r]^(.58){\theta}
\ar[d]_{\varepsilon} & \SS[\pi]\pco \ar[d]^{l}
\\
\SS[\pi]\pco \ar[r]_(.47){l} & \HH \ZZ[\pi]\pco.
 }
\end{split}
\end{equation}
\end{Le}

\begin{proof}
The maps in the diagram are defined as follows. The map $l$ is the
linearization map. The map $\varepsilon$ is induced by the
projection $SG \times \pi \ra \pi$. To obtain the map $\theta$ when
$\pi = 1$ take the adjoint of the inclusion $SG \hookrightarrow Q
(S^0)$, which is a map of the form $\SS[SG] \ra \SS$. For $\pi \neq
1$ adjoin $\pi$ to this map.

If $\pi = 1$ the diagram (\ref{square1-diagram})  is just the
diagram from \cite[Proposition 9.11]{Boekstedt-Madsen(1994)}. But
adjoining $\pi$ is the same as taking the wedge of $\pi$-many copies
of this diagram for $\pi =1$. So if the diagram for $\pi = 1$ is
$(4p-6)$-cartesian after $p$-completion then also the diagram with
adjoined $\pi$ is $(4p-6)$-cartesian after $p$-completion.
\end{proof}

\begin{proof}[Proof of Proposition \ref{K-thy-trick-BM}]
By Proposition \ref{TC-properties} applying $\TC (-,p)\pco$ to the
diagram of Lemma \ref{square1} yields a $(4p-7)$-cartesian square.
Recalling $SG \simeq \Omega BSG$ and denoting $\TC (X,p) = \TC
(\SS[\Omega X],p)$ we obtain a $(4p-7)$-cartesian square
\begin{equation} \label{TC-square}
\begin{split}
\xymatrix{ \TC(BSG \times B\pi,p)\pco \ar[r] \ar[d] & \TC(B\pi,p)\pco \ar[d] \\
\TC(B\pi,p)\pco \ar[r] & \TC (\ZZ[\pi],p)\pco.}
\end{split}
\end{equation}
We have a similar square with $K$-theory instead of $\TC (-,p)$
which maps via the cyclotomic trace map $\trc$ into the square
(\ref{TC-square}). By Theorem \ref{Dundas} the sides of the cube
containing both $K$-theory and $\TC(-,p)$ are homotopy cartesian and
hence also the whole cube is homotopy cartesian. But then the side
consisting of $K$-theories is as cartesian as the opposite side
which is the diagram (\ref{TC-square}).
\end{proof}

\begin{proof}[Proof of Proposition \ref{F-K-S^1-calculation}]
1. This follows immediately from Proposition \ref{F-K-Bpi-calc}
since $S^1 = B\ZZ$ and from the observation that the spectrum
$\Sigma^{2p-2} \HH \FF_p \wedge (S^1 (n)_+)$ already is
$p$-complete. For part 2. we need some notation. Recall that we have
a homotopy equivalence
\begin{equation} \label{Lambda-S^1}
\Lambda S^1 \simeq \coprod_{n \in \ZZ} S^1(n).
\end{equation}
By $d \alpha_1 \in \pi_{2p-2} BSG = H_{2p-2} (BSG)$, $1_n \in H_0
(S^1 (n))$, $s_n \in H_1 (S^1 (n)$ are denoted the generators of the
respective groups. Further we denote
\begin{align*} \beta_n = 1_n \cdot d \alpha_1 & \in
\pi_{2p-2} (S^1(n)_+ \wedge BSG)
\cong H_{2p-2} (S^1(n)_+ \wedge BSG) \\
\gamma_n = s_n \cdot d \alpha_1 & \in \pi_{2p-1} (S^1(n)_+ \wedge
BSG) \cong H_{2p-1} (S^1(n)_+ \wedge BSG).
\end{align*}

2. The Verschiebung operations give $\pi_\ast F_K(S^1)\pco$ a
structure of a $\ZZ_p [\NN_\times]$-module. Recall that they are
induced by maps of $\SS$-algebras $\varphi_n \colon \SS [M \times
\ZZ] \ra \SS [M \times \ZZ]$. These induce the Verschiebung
operations also on $\Sigma^\infty_+ \Lambda B(M \times \ZZ) \simeq
\Sigma^\infty_+ \Lambda X \times \Lambda S^1$, where $M = \Omega X$.
Therefore the natural map
\[
A^{\fd} (X \times S^1) \simeq K(\SS [M \times \ZZ]) \ra
\Sigma^\infty_+ \Lambda B(M \times \ZZ) \simeq \Sigma^\infty_+
\Lambda X \times \Lambda S^1
\]
is compatible with the Verschiebung operation. In more detail, the
operation $V_n$ is induced by $\varphi_n \colon \ZZ \ra \ZZ$ given
by multiplication by $n$. This induces an $n$-fold cover map on
$S^1$ (which is the target $S^1$ in $\Lambda S^1$). The effect on
$\tilde{H}_\ast (\Lambda S^1_+; \FF_p)$ is $1_m \mapsto 1_{nm}$ and
$s_m \mapsto n s_{nm}$.
\end{proof}

\begin{proof}[Proof of Corollary \ref{ZZ[NN]-module-F_K}]
This follows immediately from Proposition \ref{F-K-S^1-calculation}
(2).
\end{proof}

\begin{proof}[Proof of Proposition \ref{pi-Nil-calculation} and Corollary \ref{Wh(Bpi)}]
This follows from the homotopy fibration sequence
(\ref{lin-assembly}) and the calculation of \cite[Proposition
9.11]{Boekstedt-Madsen(1994)} that $F_K(\ast)\pco$ is
$(2p-3)$-connected and $\pi_{2p-2} F_K(\ast)\pco = \FF_p \{ d
\alpha_1 \}$. Further we have $\pi_\ast S^1_+ \wedge F_K(\ast)\pco =
\FF_p \{ \beta_0, \gamma_0 \}$ where $\beta_0 = 1_0 \cdot d
\alpha_1$, $\gamma_0 = s_0 \cdot d \alpha_1$. The same argument
works for general discrete $\pi$.
\end{proof}

\begin{proof}[Proof of Corollary \ref{ZZ[NN]-module-nil}]
This follows from Proposition \ref{pi-Nil-calculation} and Corollary
\ref{ZZ[NN]-module-F_K}.
\end{proof}

\begin{proof}[Proof of Corollary \ref{ZZ[NN*NN]-module-F_K}]
Suppose that $\pi_{2p-1} \NAPM (*,M)\pco$ is finitely generated as a
$\zz_p[\nn_\times \ast \nn_\times]$-module with the generating set
$A=\{a_1,\ldots , a_k\}$ and assume $a_1=\gamma_1$. We claim that
the following three statements hold
\begin{align}
\bigcup_{n \in \nn_{\times} \smallsetminus\{1\}} \text{Im}(\text{V}_n(\bigoplus_{i \in \nn_\times} \FF_p \{ \gamma_i\} )) & \subset \bigoplus_{\substack{i \in \nn_{\times} \\ \text{for} \,  i \neq p^l}} \FF_p \{ \gamma_i \} \label{st-(1)} \\
|\bigcup_{\substack{n \in \nn_\times \\ a_i \in A}} \text{F}_n(a_i) \cap \bigoplus_{l \in \nn_\times} \FF_p \{\gamma_{p^l}\}| & < \infty  \label{st-(2)} \\
\bigcup_{n \in \nn_{\times} \smallsetminus\{1\}}
\text{Im}(\text{F}_n(\bigoplus_{\substack{i \in \nn_\times \\
\text{for} \,  i \neq  p^l}} \FF_p \{ \gamma_i\} )) & \subset
\bigoplus_{\substack{i \in \nn_\times \\ \text{for} \,  i \neq p^l}}
\FF_p \{ \gamma_i \} \cup \bigcup_{n \in \nn_\times} \text{F}_n
(\FF_p \{\gamma_1\}). \label{st-(3)}
\end{align}
The statement (\ref{st-(1)}) follows immediately from Proposition
\ref{F-K-S^1-calculation} (2). The statement (\ref{st-(2)}) follows
since for  any $x$ we have $F_n(x) = 0$ when $n \geq N(x)$ for some
$N(x)$. To prove (\ref{st-(3)}) let
\[
\gamma_m \in \bigoplus_{\substack{i \in \nn \\ \text{for} \, i \neq
p^l}} \FF_p \{ \gamma_i\}.
\]
Again by Proposition \ref{F-K-S^1-calculation} (2) we have
$\gamma_m= \text{V}_m(z \gamma_1)$ for some $z \in \FF_p$. Then
\begin{align*}
\text{F}_n (\gamma_m ) & = \text{F}_n \text{V}_m (z \gamma_1 ) \\
                         & = \text{gcd}(m,n) \text{F}_{n'} \text{V}_{m'} (z \gamma_1) \\
                         & = \text{gcd}(m,n) \text{V}_{m'} \text{F}_{n'} (z \gamma_1 )
\end{align*}
where $\text{gcd}(m,n)$ denotes the greatest common divisor, $m =
\text{gcd}(m,n) \cdot m'$ and $n = \text{gcd}(m,n) \cdot n'$. If $m'
=1$ then $\text{F}_n (\gamma_m ) \in \bigcup_{n \in \nn_\times}
\text{F}_n (\FF_p \{\gamma_1\})$ and if $m' > 1$ then by
(\ref{st-(1)}) we have $\text{F}_n (\gamma_m ) \in
\bigoplus_{\substack{i \in \nn_\times \\ \text{for} \,  i \neq p^l}}
\FF_p \{ \gamma_i \}$.

By (\ref{st-(1)}) and (\ref{st-(3)}) the subgroup
\[
\bigoplus_{\substack{i \in \nn_{\times} \\ \text{for} \,  i \neq
p^l}} \FF_p \{ \gamma_i \} \oplus \bigcup_{\substack{n \in
\nn_\times \\ a_i \in A}} \text{F}_n(a_i)
\]
is invariant under the Frobenius and Verschiebung operations.
Furthermore all the generators are contained in this subgroup. But
this gives a contradiction since by (\ref{st-(2)}) there are
elements in the complement of this subgroup in $\oplus_{n \in
\nn_\times} \FF_p \{\gamma_n\}$.
\end{proof}


\small
\bibliography{operations}

\newcommand{\etalchar}[1]{$^{#1}$}
\begin{thebibliography}{HKV{\etalchar{+}}02}

\bibitem[BM94]{Boekstedt-Madsen(1994)}
M.~B{\"o}kstedt and I.~Madsen.
\newblock Topological cyclic homology of the integers.
\newblock {\em Ast\'erisque}, (226):7--8, 57--143, 1994.
\newblock $K$-theory (Strasbourg, 1992).

\bibitem[CdS95]{Conolly-et-al(1995)}
Francis~X. Connolly and M{\'a}rio O.~M. da~Silva.
\newblock The groups {$N\sp rK\sb 0(\mathbf Z\pi)$} are finitely generated
  {$\mathbf Z[\mathbf N\sp r]$}-modules if {$\pi$} is a finite group.
\newblock {\em $K$-Theory}, 9(1):1--11, 1995.

\bibitem[DGM04]{Dundas-Goodwillie-McCarthy(2004)}
B.~Dundas, T.~Goodwillie, and R.~McCarthy.
\newblock {\em The local structure of the algebraic K-theory}.
\newblock Preprint, 2004.

\bibitem[Dun97]{Dundas(1997)}
Bj{\o}rn~Ian Dundas.
\newblock Relative {$K$}-theory and topological cyclic homology.
\newblock {\em Acta Math.}, 179(2):223--242, 1997.

\bibitem[Far77]{Farrell}
F.~T. Farrell.
\newblock The nonfiniteness of {N}il.
\newblock {\em Proc. Amer. Math. Soc.}, 65(2):215--216, 1977.

\bibitem[FJ91]{Farrell-Jones(1991)}
F.~T. Farrell and L.~E. Jones.
\newblock Stable pseudoisotopy spaces of compact non-positively curved
  manifolds.
\newblock {\em J. Differential Geom.}, 34(3):769--834, 1991.

\bibitem[Goo90]{Goodwillie(1990)}
Thomas~G. Goodwillie.
\newblock Calculus. {I}. {T}he first derivative of pseudoisotopy theory.
\newblock {\em $K$-Theory}, 4(1):1--27, 1990.

\bibitem[Gra76]{Grayson(1976)}
Daniel Grayson.
\newblock Higher algebraic {$K$}-theory. {II} (after {D}aniel {Q}uillen).
\newblock In {\em Algebraic $K$-theory (Proc. Conf., Northwestern Univ.,
  Evanston, Ill., 1976)}, pages 217--240. Lecture Notes in Math., Vol. 551.
  Springer, Berlin, 1976.

\bibitem[HKV{\etalchar{+}}01]{fund1}
T.~H{\"u}ttemann, J.~Klein, W.~Vogell, F.~Waldhausen, and B.~Williams.
\newblock The ``fundamental theorem'' for the algebraic {$K$}-theory of spaces.
  {I}.
\newblock {\em J. Pure Appl. Algebra}, 160(1):21--52, 2001.

\bibitem[HKV{\etalchar{+}}02]{fund2}
T.~H{\"u}ttemann, J.~Klein, W.~Vogell, F.~Waldhausen, and B.~Williams.
\newblock The ``fundamental theorem'' for the algebraic {$K$}-theory of spaces.
  {II}. {T}he canonical involution.
\newblock {\em J. Pure Appl. Algebra}, 167(1):53--82, 2002.

\bibitem[Igu82]{Igusa(1982)}
Kiyoshi Igusa.
\newblock On the algebraic {$K$}-theory of {$A\sb{\infty }$}-ring spaces.
\newblock In {\em Algebraic $K$-theory, Part II (Oberwolfach, 1980)}, volume
  967 of {\em Lecture Notes in Math.}, pages 146--194. Springer, Berlin, 1982.

\bibitem[KR97]{Klein-Rognes(1997)}
John~R. Klein and John Rognes.
\newblock The fiber of the linearization map {$A(*)\to K({\bf Z})$}.
\newblock {\em Topology}, 36(4):829--848, 1997.

\bibitem[KW95]{fund3}
J.~Klein and B.~Williams.
\newblock The ``fundamental theorem'' for the algebraic {$K$}-theory of spaces.
  {III}. {T}he nil-term.
\newblock {\em Preprint}, 1995.

\bibitem[Lod81]{Loday(1980)}
Jean-Louis Loday.
\newblock On the boundary map {$K\sb{3}(\Lambda /I)\rightarrow K\sb{2}(\Lambda
  ,\,I)$}.
\newblock In {\em Algebraic $K$-theory, Evanston 1980 (Proc. Conf.,
  Northwestern Univ., Evanston, Ill., 1980)}, volume 854 of {\em Lecture Notes
  in Math.}, pages 262--268. Springer, Berlin, 1981.

\bibitem[Lyd99]{Lydakis(1999)}
Manos Lydakis.
\newblock Smash products and {$\Gamma$}-spaces.
\newblock {\em Math. Proc. Cambridge Philos. Soc.}, 126(2):311--328, 1999.

\bibitem[Mad94]{Madsen(1995)}
Ib~Madsen.
\newblock Algebraic {$K$}-theory and traces.
\newblock In {\em Current developments in mathematics, 1995 (Cambridge, MA)},
  pages 191--321. Int. Press, Cambridge, MA, 1994.

\bibitem[Rog02]{Rognes(2002)}
John Rognes.
\newblock Two-primary algebraic {$K$}-theory of pointed spaces.
\newblock {\em Topology}, 41(5):873--926, 2002.

\bibitem[Rog03]{Rognes(2003)}
John Rognes.
\newblock The smooth {W}hitehead spectrum of a point at odd regular primes.
\newblock {\em Geom. Topol.}, 7:155--184 (electronic), 2003.

\bibitem[Ros94]{Rosenberg(1994)}
Jonathan Rosenberg.
\newblock {\em Algebraic {$K$}-theory and its applications}, volume 147 of {\em
  Graduate Texts in Mathematics}.
\newblock Springer-Verlag, New York, 1994.

\bibitem[Seg74]{Segal}
G.~Segal.
\newblock Categories and cohomology theories.
\newblock {\em Topology}, 13:293--312, 1974.

\bibitem[Tod58]{Toda(1958)}
Hirosi Toda.
\newblock {$p$}-primary components of homotopy groups. {I}. {E}xact sequences
  in {S}teenrod algebra.
\newblock {\em Mem. Coll. Sci. Univ. Kyoto. Ser. A. Math.}, 31:129--142, 1958.

\bibitem[Wal78]{Waldhausen(1978)}
Friedhelm Waldhausen.
\newblock Algebraic {$K$}-theory of topological spaces. {I}.
\newblock In {\em Algebraic and geometric topology (Proc. Sympos. Pure Math.,
  Stanford Univ., Stanford, Calif., 1976), Part 1}, Proc. Sympos. Pure Math.,
  XXXII, pages 35--60. Amer. Math. Soc., Providence, R.I., 1978.

\bibitem[Wal85]{Waldhausen1126}
F.~Waldhausen.
\newblock Algebraic {$K$}-theory of spaces.
\newblock In {\em Algebraic and geometric topology (New Brunswick, N.J.,
  1983)}, volume 1126 of {\em Lecture Notes in Math.}, pages 318--419.
  Springer, Berlin, 1985.

\bibitem[WW01]{Weiss-Williams(2001)}
Michael Weiss and Bruce Williams.
\newblock Automorphisms of manifolds.
\newblock In {\em Surveys on surgery theory, Vol. 2}, volume 149 of {\em Ann.
  of Math. Stud.}, pages 165--220. Princeton Univ. Press, Princeton, NJ, 2001.

\end{thebibliography}
\bibliographystyle{alpha}

\end{document}